\documentclass[a4paper,11pt]{amsart}
\usepackage[french]{babel}

\usepackage{amscd,amsfonts,amsmath,amssymb}

\input xy
\xyoption{all}
\usepackage[all]{xy}

\begin{document}
\title{Sur une conjecture de Mukai}
\author{Laurent BONAVERO, Cinzia CASAGRANDE, Olivier DEBARRE, St\'ephane
DRUEL}
\date{2002}
\maketitle

\def\restriction{\string |}
\def\cad{c'est-\`a-dire}
\def\emptyset{\varnothing}
\def\moins{\mathop{\hbox{\vrule height 3pt depth -2pt
width 5pt}\,}}
\def\setminus{\moins}

\newcommand{\pp}{\rm ppcm}
\newcommand{\pg}{\rm pgcd}
\newcommand{\Ker}{\rm Ker}
\newcommand{\C}{{\mathbb C}}
\newcommand{\Q}{{\mathbb Q}}
\newcommand{\GL}{\rm GL}
\newcommand{\SL}{\rm SL}
\newcommand{\diag}{\rm diag}

\newcommand{\ev}{\operatorname{ev}}
\def\refname{R\'ef\'erences}
\def\finpreuve
{\hskip 3pt \vrule height6pt width6pt depth 0pt}

\newtheorem{theo}{Th\'eor\`eme}
\newtheorem{prop}{Proposition}
\newtheorem{lemm}{Lemme}
\newtheorem{lemmf}{Lemme fondamental}
\newtheorem{defi}{D\'efinition}
\newtheorem{exo}{Exercice}
\newtheorem{rem}{Remarque}
\newtheorem{cor}{Corollaire}
\newcommand{\CC}{{\mathbb C}}
\newcommand{\ZZ}{{\mathbb Z}}
\newcommand{\RR}{{\mathbb R}}
\newcommand{\QQ}{{\mathbb Q}}
\newcommand{\FF}{{\mathbb F}}
\newcommand{\PP}{{\mathbb P}}
\newcommand{\codim}{\operatorname{codim}}
\newcommand{\Exc}{\operatorname{Exc}}
\newcommand{\Ho}{\operatorname{Hom}}
\newcommand{\Pic}{\operatorname{Pic}}
\newcommand{\Aut}{\operatorname{Aut}}
\newcommand{\NE}{\operatorname{NE}}
\newcommand{\Nun}{\operatorname{N}}
\newcommand{\card}{\operatorname{card}}
\newcommand{\Hilb}{\operatorname{Hilb}}
\newcommand{\mult}{\operatorname{mult}}
\newcommand{\vol}{\operatorname{vol}}
\newcommand{\divi}{\operatorname{div}}
\newcommand{\pr}{\operatorname{pr}}
\newcommand{\con}{\operatorname{cont}}
\newcommand{\Spec}{\operatorname{Spec}}
\newcommand{\ima}{\operatorname{Im}}
\newcommand{\Chow}{\operatorname{Chow}}
\newcommand{\lieu}{\operatorname{lieu}}
\newcommand{\ol}[1]{\ensuremath{\mathcal{O}_{#1}}}
\newcommand{\Div}{\ensuremath{\operatorname{Div}}}
\newcommand{\Supp}{\operatorname{Supp}}
\newcommand{\Star}{\operatorname{Star}}
\newcommand{\ph}{\ensuremath{\varphi}}
\newcommand{\lis}[1]{\ensuremath{\{x_1,\dots,x_{#1}}\}}
\newcommand{\som}[1]{\ensuremath{x_1+\dots+x_{#1}}}
\newcommand{\fx}{\ensuremath{\Sigma_X}}
\newcommand{\fy}{\ensuremath{\Sigma_Y}}
\newcommand{\Z}{\ensuremath{\mathbb{Z}}}
\newcommand{\f}{\ensuremath{\Sigma}}
\newcommand{\Int}{\operatorname{Int}}
\newcommand{\N}{\ensuremath{\mathcal{N}_1}}
\newcommand{\A}{\ensuremath{\mathcal{A}_1}}
\newcommand{\PC}{\operatorname{PC}}

\newtheorem{teo}{Th\'eor\`eme}
\newtheorem{lemma}[teo]{Lemme}
\newtheorem{conj}[teo]{Conjecture}

\newcounter{subsub}[subsection]

\def\thesubsub{\thesubsection .\arabic{subsub}}
\def\subsub#1{\addtocounter{subsub}{1}\par\vspace{3mm}
\noindent{\bf \thesubsub ~ #1 }\par\vspace{2mm}}
\def\coker{\mathop{\rm coker}\nolimits}
\def\pr{\mathop{\rm pr}\nolimits}
\def\im{\mathop{\rm Im}\nolimits}
\def\hfl#1#2{\smash{\mathop{\hbox to 12mm{\rightarrowfill}}
\limits^{\scriptstyle#1}_{\scriptstyle#2}}}
\def\vfl#1#2{\llap{$\scriptstyle #1$}\big\downarrow
\big\uparrow
\rlap{$\scriptstyle #2$}}
\def\diagram#1{\def\normalbaselines{\baselineskip=0pt
\lineskip=10pt\lineskiplimit=1pt}   \matrix{#1}}
\def\limind{\mathop{\oalign{lim\cr
\hidewidth$\longrightarrow$\hidewidth\cr}}}

\long\def\InsertFig#1 #2 #3 #4\EndFig{
\hbox{\hskip #1 mm$\vbox to #2 mm{\vfil\includegraphics{#3}}#4$}}
\long\def\LabelTeX#1 #2 #3\ELTX{\rlap{\kern#1mm\raise#2mm\hbox{#3}}}

\newcommand{\num}{\stepcounter{subsection}{\bf\thesubsection.}}
\newcommand{\marque}{\stepcounter{subsection}{\thesubsection. }}
\newcommand{\hs}{\hspace*{5mm}}
\newcommand{\voir}{{\it cf.\ }}

{\let\thefootnote\relax
\footnote{
\textrm{Mots cl\'es: vari\'et\'es de Fano,
th\'eorie de Mori, g\'eom\'etrie torique}.\\
\hspace*{4.4mm}\textrm{Classification~A.M.S.: 14J45, 14E30,
14M25.}
}}

\begin{center}
\begin{minipage}{130mm}
\scriptsize

{\bf R\'esum\'e.} G\'en\'eralisant une question de Mukai,
nous conjecturons qu'une vari\'et\'e
de Fano $X$ de nombre de Picard $\rho_X$ et
de pseudo-indice $\iota_X$ v\'erifie $\rho_X (\iota_X-1) \le \dim(X)$.
Nous d\'emontrons cette conjecture dans plusieurs situations:
$X$ est une vari\'et\'e
de Fano de dimension $\le4$,
$X$ est une vari\'et\'e
de Fano torique  de dimension $\le7$
ou $X$ est une vari\'et\'e
de Fano toriques de dimension arbitraire avec $\iota_X \ge \dim(X) /3+1$.
Enfin, nous pr\'esentons une approche nouvelle pour le cas g\'en\'eral.

\medskip

{\bf Abstract.} Generalizing a question of Mukai,
we conjecture that a Fano manifold $X$
with Picard number $\rho_X$ and pseudo-index $\iota_X$
satisfies $\rho_X (\iota_X-1) \le \dim(X)$.
We prove this inequality in several situations:
$X$ is a   Fano manifold  of dimension $\le4$, $X$ is a toric Fano
manifold  of dimension $\le7$ or $X$ is a toric Fano manifold of
arbitrary dimension with $\iota_X \ge \dim(X) /3+1$. Finally, we offer a
new approach to the general case.

\end{minipage}
\end{center}

\section*{Introduction}

\hs Soit $X$ une vari\'et\'e de Fano, \cad\
une vari\'et\'e projective lisse dont le fibr\'e anticanonique
$-K_X  $ est ample, d\'efinie sur $\CC$.
Notons $\iota_X$ le pseudo-indice de $X$, c'est-\`a-dire
le plus petit entier de la forme $-K_X \cdot C$ o\`u $C$ est une courbe
rationnelle de $X$, et $\rho _X$ le nombre de Picard
de $X$. L'indice $r_X$ de   $X$
est le plus grand entier $m$ tel qu'il existe un fibr\'e en droites
$L$ satisfaisant $-K_X = mL$ dans $\Pic (X)$; c'est un nombre
particuli\`erement adapt\'e
\`a l'\'etude des vari\'et\'es de Fano $X$ avec $\rho _X =1$, tandis
le pseudo-indice nous semble plus adapt\'e
\`a l'\'etude des vari\'et\'es de Fano $X$ avec $\rho _X \ge 2$
(penser par exemple aux produits d'espaces projectifs).
Mukai a propos\'e dans \cite{Muk88}
l'in\'egalit\'e $\rho_X (r_X-1) \le \dim(X)$,
v\'erifi\'ee  en dimension $\le3$ \`a l'aide de la classification.
Nous proposons plus g\'en\'eralement l'in\'egalit\'e
$$ \quad\quad\rho_X (\iota_X-1) \le \dim(X)\leqno (*)$$
avec \'egalit\'e si et seulement si $X\simeq
(\PP^{\iota_{X}-1})^{\rho_X}$, en partie motiv\'es par le

\medskip

\noindent {\bf Th\'eor\`eme (\cite{wisn2}).} {\em
Soit $X$ une vari\'et\'e de Fano de dimension $n$,
d'indice $r_X$ et de pseudo-indice $\iota_X$.
Si $ 2 \iota_X > n +2$, on a $\rho_X = 1$.
Si $2 r_X = n+2$, on a $\rho_X = 1$ sauf si $X\simeq
(\PP^{\iota_{X}-1})^2$. }

\medskip

\hs D'apr\`es la th\'eorie de Mori, une vari\'et\'e de Fano $X$ de
dimension $n$ v\'erifie toujours $\iota_X \le n+1$. Si $\iota_X = n+1$,
on a $X \simeq \PP ^n$ d'apr\`es un r\'esultat de Cho, Miyaoka et
Shepherd-Barron (\cite{CMS00}) dont une d\'emonstration simplifi\'ee se
trouve dans
\cite{Keb01}.

\hs Dans ce travail, on se propose de montrer
l'in\'egalit\'e ($\ast$) lorsque $X$ est une vari\'et\'e de Fano
de dimension $3$ ou $4$ (l'in\'egalit\'e ($\ast$) est
v\'erifi\'ee
directement en dimension $2$)
et lorsque $X$ est une vari\'et\'e de Fano torique
de pseudo-indice au moins \'egal
\`a $\dim(X) /3 +1$.\\
\hs L'in\'egalit\'e ($\ast$) peut \^etre v\'erifi\'ee directement pour
les vari\'et\'es de Fano de dimension $4$ et d'indice au moins $2$,
dont la classification est \'etablie
par exemple dans \cite{IsP99}.

\hs Nous d\'emontrons, sans utiliser les r\'esultats de classification
des vari\'et\'es de Fano de dimension $3$ ou $4$, le r\'esultat suivant,
nouveau dans le cas des vari\'et\'es d'indice
$1$ et de pseudo-indice au moins $2$:

\medskip

\noindent {\bf Th\'eor\`eme.} {\em
Si $X$ est une vari\'et\'e de Fano de dimension $3$ ou $4$,
de pseudo-indice $\iota_X$ et de nombre de Picard $\rho_X$, on a
$\rho_X (\iota_X-1) \le\dim(X)$. Si $\rho_X (\iota_X-1)=\dim(X)$, on a
 $X\simeq(\PP^{\iota_{X}-1})^{\rho_{X}}$.
}

\medskip

\hs Dans la situation torique, nos r\'esultats principaux sont:

\medskip

\noindent {\bf Th\'eor\`eme.} {\em
Soit $X$ une vari\'et\'e de Fano torique de dimension $n$,
de pseudo-indice $\iota_X$ et de nombre de Picard $\rho_X$.
Si $3 \iota_X\ge n+3$, on a $\rho_X(\iota_X-1)\le n$. Si de plus
$\rho_X(\iota_X-1)=n$, on a
$X\simeq (\PP ^{\iota_X-1})^{\rho_X}$.
}

\medskip

\noindent {\bf Th\'eor\`eme.} {\em
Soit $X$ une vari\'et\'e de Fano torique de dimension $n \le 7$,
de pseudo-indice $\iota_X$ et de nombre de Picard $\rho_X$.
On a $\rho_X(\iota_X-1)\le n$. Si
$\rho_X(\iota_X-1)=n$, on a
$X\simeq (\PP ^{\iota_X-1})^{\rho_X}$.
}

\medskip

\hs Les m\'ethodes utilis\'ees sont celles de la th\'eorie de Mori,
dont nous rappelons plus bas quelques uns des
r\'esultats ou d\'efinitions.
Mentionnons les deux r\'esultats suivants obtenus au cours de
la d\'emonstration des th\'eor\`emes pr\'ec\'edents et qui ont
leur int\'er\^et propre.

\medskip

\noindent {\bf Proposition.} {\em
Soit $X$ une vari\'et\'e de Fano de dimension $n$,
de nombre de Picard $\rho_{X}\ge n$ et de pseudo-indice $\iota_{X} \ge
2$. Si toutes les contractions extr\'emales de $X$  sont des fibrations,
$X$ est isomorphe \`a $(\mathbb{P}^{1})^{n}$. }

\medskip

\noindent {\bf Th\'eor\`eme.} {\em
Soit $X$ une vari\'et\'e
de Fano de dimension $n\ge 4$.
Si toutes les contractions extr\'emales de $X$
sont ou bien des \'eclatements lisses de centre une courbe
lisse ou bien des fibrations
de dimension relative $1$ et que $X$ poss\`ede au moins
une contraction extr\'emale birationnelle, on a  $\rho_X \le 3$.
}

\medskip

\hs Nous pr\'esentons enfin une approche
possible, en toute dimension, o\`u l'in\'egalit\'e
($\ast$) d\'ecoulerait de l'existence de certaines cha\^ines
de courbes rationnelles. Le r\'esultat suivant fournit une motivation
suppl\'ementaire en direction de la conjecture de Mukai g\'en\'eralis\'ee.

\medskip

\noindent {\bf Th\'eor\`eme.} {\em
Soit $X$ une vari\'et\'e de Fano de dimension
$n$, de nombre de Picard $\rho_X$ et de pseudo-indice $\iota_X$. S'il existe
des familles propres irr\'e\-ductibles $V^1,\dots,V^{\rho_X}$
de cour\-bes rationnelles irr\'e\-ductibles sur $X$ dont les
classes dans $\Nun_{1}(X)_{\QQ}$ sont lin\'eairement
ind\'ependantes et des courbes $C^1,\dots,C^{\rho_{X}}$
avec $[C^j] \in V^j$ pour tout $1\le j\le \rho_{X}$
et $C^j \cap C^{j+1} \neq \emptyset$ pour tout $1\le j\le \rho_{X}-1$,
on a
$\rho_X (\iota_X -1)\le n$.
}

\section{Notations}
\marque Si $X$ est une vari\'et\'e projective, on note
$$\Nun_1(X)_{\QQ} = \{ \sum_i a_i C_i \mid a_i \in \QQ  ,
C_i  \ \mbox{courbe irr\'eductible de}  \ X \}/\equiv $$
o\`u $\equiv$ d\'esigne l'\'equivalence num\'erique.
Le c\^one de Mori, ou c\^one des courbes effectives,
est le sous-c\^one de $\Nun_1(X)_{\QQ}$ d\'efini par
$$ \NE (X) = \{ Z\in \Nun_1(X)_{\QQ} \mid Z \equiv
\sum_i a_i C_i  ,  a_i\ge 0\}.$$
Si $E$ est un diviseur  et $R$ une ar\^ete
de $\NE(X)$, la notation $E\cdot R >0$ (resp. $E\cdot R =0$)
signifie
$E \cdot C >0$ (resp. $E\cdot C =0$)
pour toute courbe $C$ de $X$ telle que $[C] \in R$.
La partie de $\NE (X)$ situ\'ee dans le demi-espace
ouvert $\{ -K_X >0\}$ est localement po\-ly\'e\-dra\-le
et pour
toute ar\^ete $R$ de $\NE (X)$ telle que $-K_X \cdot R >0$, il existe
un morphisme $\varphi_R : X \to X_R$ \`a fibres connexes, appel\'e
contraction extr\'emale,
de $X$ sur une vari\'et\'e
projective normale $X_R$ tel que
les courbes irr\'eductibles contract\'ees par
$\varphi_R$ sont exactement celles dont la classe dans
$\Nun_1(X)_{\QQ}$ appartient \`a $R$.
\begin{itemize}
\item[$\bullet$] Soit $\dim(X_R)<\dim(X)$; on dit que $\varphi_R$ est
une fibration, et $\rho_{X_R}=\rho_X-1$.
\item[$\bullet$] Soit $\varphi_R$ est birationnelle, son lieu
exceptionnel est un diviseur irr\'eductible $E$ tel que $E\cdot R<0$:
on dit que $\varphi_R$ est une contraction divisorielle, et
$\rho_{X_R}=\rho_X-1$.
\item[$\bullet$] Soit $\varphi_R$ est birationnelle, son lieu
exceptionnel est de codimension $>1$ dans $X$:
on dit que $\varphi_R$ est une petite contraction.
\end{itemize}
Si $X$ est une vari\'et\'e de Fano, le c\^one $\NE (X)$ est
poly\'edral et si $\varphi_R : X \to X_R$ est une contraction extr\'emale,
$\NE (X_R)$ est engendr\'e par les $\varphi_R (R')$ o\`u
$R'$ d\'ecrit l'ensemble (fini) des ar\^etes de $\NE (X)$.

\medskip

\marque Rappelons un r\'esultat de
Wi\'sniewski (\cite{Wis91}, Theorem (1.1)) sur le lieu d'une contraction
extr\'emale. Soit $\varphi_R : X \to X_R$ une telle contraction,
d'ar\^ete $R$. Notons
$E_{R}$ le lieu de $X$ couvert par les courbes
contract\'ees par $\varphi_R$ et
$$ l(R) = \min \{ -K_X\cdot C \mid C \
\mbox{courbe rationnelle contract\'ee par }
\varphi_R \}$$ la longueur de $R$. Pour toute composante
irr\'eductible $F$
d'une fibre non triviale de $\varphi_R$, on a
$$ 2\dim (E_{R})\ge\dim (E_{R}) + \dim (F) \ge \dim (X) +
l(R) - 1\ge\dim(X)+\iota_X-1.$$

\marque Nous utiliserons aussi
le r\'esultat suivant: soient $X$ une vari\'et\'e de Fano, $R$ une ar\^ete
de $\NE (X)$ et $\varphi_R : X \to X_R$ la contraction associ\'ee.
{\em Si toutes les fibres de $\varphi_R $ sont de dimension $1$},
 $X_R$ est lisse (\cite{And85}, Theorem
3.1(ii)) et une fibre g\'en\'erale $F$ de
$\varphi_R$ est une courbe rationnelle lisse qui v\'erifie $-K_X\cdot
F=2\ge\iota_X$. Deux cas sont alors
possibles:
\begin{itemize}
\item[$\bullet$]soit $\varphi_R$ est une  fibration en $\PP ^1$ et
$X_R$ est une   vari\'et\'e de Fano (\cite{kmm}, Corollary 2.9);
\item[$\bullet$]soit $\varphi_R$ a au moins une fibre singuli\`ere
 et $\iota_X=1$.
\end{itemize}

\section{La dimension trois}
\marque V\'erifions l'in\'egalit\'e ($\ast$) en dimension $3$.
Le seul cas non trivial et non couvert
par le r\'esultat de Wi\'sniewski est donn\'e par la

\medskip

\noindent {\bf Proposition \num} {\em
Soit $X$ une vari\'et\'e de Fano de dimension $3$.
Si $\iota_X =2$, on a $\rho_X \le 3$.
}

\medskip

\noindent {\em D\'emonstration.} Supposons $\rho_X \ge 4$. La vari\'et\'e $X$
poss\`ede  au moins $4$ contractions extr\'emales distinctes qui
ne peuvent pas toutes \^etre des fibrations (\cite{Wis91},
Theorem (2.2)).  Les contractions
extr\'emales sont compl\`etement d\'ecrites par Mori (\cite{Mor82},
Theorem 3.3) et
seule la contraction birationnelle lisse de centre un point
est de longueur au moins $2$. Cette possibilit\'e est par exemple exclue
par la classification des vari\'et\'es dont
l'\'eclat\'ee en un point est de Fano (\cite{BCW01}).\finpreuve

\medskip

\hs  Toujours en dimension $3$, le cas d'\'egalit\'e
dans l'in\'egalit\'e ($\ast$) est caract\'eris\'e par la

\medskip

\noindent {\bf Proposition \num} {\em
Si $\rho_X(\iota_{X}-1)=\dim (X)=3$, on a
$X\simeq (\mathbb{P}^{\iota_{X}-1})^{\rho_{X}}.$ }

\medskip

\noindent {\em D\'emonstration.} Si $\rho_{X}=1$
et $\iota_{X}=4$, la vari\'et\'e $X$ est isomorphe \`a $ \PP^{3}$
(\cite{CMS00}, \cite{Keb01}). Supposons
$\rho_{X}=3$ et $\iota_{X}=2$ et consid\'erons une contraction extr\'emale
$\varphi:X\to Z$. Si $\varphi$
est divisorielle de lieu exceptionnel $E$, alors
$\varphi(E)$ est une courbe (\voir 1.2) et $\varphi$
s'identifie \`a l'\'eclatement d'une courbe lisse
dans la vari\'et\'e lisse $Z$ (\cite{Mor82},
Corollary 3.4.1), ce qui est absurde
puisque $\iota_{X}=2$. Les contractions extr\'emales sont donc
toutes des fibrations et le r\'esultat
annonc\'e est donn\'e par la proposition qui suit.\finpreuve

\medskip

\noindent {\bf Proposition \num} {\em
Soit $X$ une vari\'et\'e de Fano de dimension $n$,
de nombre de Picard $\ge n$ et de pseudo-indice    $\ge2$. Si toutes les
contractions extr\'emales de $X$ sont des fibrations, $X$ est
isomorphe
\`a
$(\mathbb{P}^{1})^{n}$.}

\medskip

{\em D\'emonstration.} Montrons la proposition par
r\'ecurrence sur  $n$,
le cas $n=1$ \'etant imm\'ediat. Par \cite{Wis91}, Theorem (2.2),
  toutes les contractions extr\'emales sont de
dimension relative $1$ et  il y en
a exactement $n$: le c\^one   $\NE(X)$ est donc simplicial.
Notons $R_1,\dots,R_n$ ses ar\^etes. Pour $1\le i\le n-1$, le c\^one
$V_{i}=R_{1}+\dots+R_{i}$ est une
face extr\'emale de $\NE(X)$ dont on note
$\varphi_i:X\to Z_{i}$ la contraction.
Par le lemme de
rigidit\'e, $\varphi_{i+1}$ se factorise en
$ X\stackrel{\varphi_i}{\to} Z_{i}
\to Z_{i+1}
$. L'inclusion $V_{i}\subset
V_{i+1}$
\'etant stricte, on a $\dim(Z_{i})>\dim(Z_{i+1})$. On en
d\'eduit
$\dim(Z_{n-1})=1$ et $Z_{n-1}$ est isomorphe
\`a $\mathbb{P}^{1}$.

\hs Pour $1\le i\le n$, on consid\`ere la contraction
$\ph_{R_i}:X\to W_i$. Les fibres de  $\varphi_{n-1}$ et de $\ph_{R_n}$ se
coupent en un nombre fini de points. Le morphisme $\ph_{R_n}$ est donc
\'equidimensionnel de dimension relative $1$.
Il en r\'esulte que $W_{n}$ est une vari\'et\'e de Fano lisse
(\voir 1.3) et
que $\ph_{R_n}$ est une fibration
  en $\PP^1$ puisque $\iota_{X} \geq 2$.
Comme toute ar\^ete de $\NE(W_{n})$ est image
d'une ar\^ete
de $\NE(X)$, le c\^one $\NE(W_{n})$ est simplicial et
toute contraction
extr\'emale de $W_{n}$ est une fibration.

\medskip

\noindent {\bf Lemme \num} {\em
Soit  $\pi: X\to Y$ une
fibration  en
$\PP^r$ entre vari\'et\'es projectives et  lisses. On suppose que
$X$ est  une vari\'et\'e de Fano. Alors,
\begin{itemize}
\item[(a)] $Y$ est
une vari\'et\'e de Fano et
$\iota_Y\ge\iota_X$;
\item[(b)]si $Y$ v\'erifie $(*)$, il en est de m\^eme
de
$X$;
\item[(c)]si $\iota_Y=\iota_X$ et que $\PP^1\to Y$ est une
courbe   de degr\'e anticanonique $\iota_Y$, le
produit $\PP^1\times_YX$  est isomorphe \`a
$\PP^1\times\PP^r$.
 \end{itemize}}
\medskip

\noindent {\em D\'emonstration.} Par \cite{kmm}, Corollary
2.9,  $Y$ est
une vari\'et\'e de Fano. Soit $ \PP^1\to Y$ une courbe
  de  degr\'e
anticanonique $\iota_Y$. Comme $H^2(\PP^1,{\mathcal O}^*_{\PP^1})$
est nul, il r\'esulte de
\cite{gro},
\S\ 8, que la fibration
$ \PP^1\times_YX\to\PP^1$ est isomorphe au projectifi\'e d'un
fibr\'e vectoriel
  $ E=
\bigoplus_{i=0}^r {\mathcal O}_{\PP^1}(a_i)$ sur $\PP^1$, o\`u les
$a_i$ sont des entiers positifs avec $a_0=0$. Si
$C_i $ est la   section de $\PP(E)\to\PP^1$  d\'efinie par le
fibr\'e quotient
${\mathcal O}_{\PP^1}(a_i)$ de $E$, on a
 $-K_{\PP(E)/\PP^1}\cdot C_i =  -ra_i  $, d'o\`u, en
notant $g$ la compos\'ee $\PP^1\to C_i\to X$,
$$\iota_X\le -K_X\cdot g_*C_i=-\pi^*K_Y\cdot g_*C_i-K_{X/Y}\cdot
g_*C_i=\iota_Y-g^*K_{X/Y}\cdot
 C_i=\iota_Y- K_{\PP(E)/\PP^1}\cdot
 C_i=
\iota_Y-ra_i
$$
ce qui prouve (a) et (c).

\hs Montrons (b). On a $\rho_Y=\rho_X-1$, d'o\`u, en utilisant
(a),
$$(\rho_X-1)(\iota_X-1)\le
\rho_Y(\iota_Y-1)\le\dim (Y)=\dim X-r.$$
Si $C$ est une
droite
 contenue dans une fibre de $\pi$, on a $-K_X\cdot
C=r+1$, de sorte que
 $\iota_X\le r+1$ et
$$\rho_X(\iota_X-1)\le \iota_X-1+ \dim
(X)-r \le \dim (X),$$
ce qui prouve (b).\finpreuve

\hs On a donc $\iota_{W_{n}} \geq 2$ et par hypoth\`ese de r\'ecurrence,
$W_{n}$ est isomorphe
\`a $(\mathbb{P}^{1})^{n-1}$. De plus, l'image r\'eciproque de toute
droite du type
$\PP^1\times\{t_2,\dots,t_{n-1}\}$ est isomorphe \`a
$(\mathbb{P}^{1})^{2}$. Puisque $\NE(X)$ est simplicial, la contraction
$\ph_{R_1+R_n}$ est la compos\'ee  de $ \ph_{R_n}$ et d'une projection
$p:W_n\to (\PP^1)^{n-2}$; de m\^eme, $\ph_{n-1}$
est la compos\'ee  de $ \ph_{R_n}$ et d'une projection
$ W_n\to  \PP^1 $. On a ainsi un diagramme commutatif
$$
\xymatrix{
&X\ar[dl]_{\varphi_{n-1}}\ar[dr]^{\ph_{R_1+R_n}}
\ar[d]^{\ph_{R_1}}
\ar[r]^-{\ph_{R_n}} &\  \ W_n\simeq(\mathbb{P}^1)^{n-1}
\ar[d]^{p} \\
\mathbb{P}^1 &\ W_1\simeq\PP^1\times(\mathbb{P}^1)^{n-2}
\ar[r]^-{p_2}\ar[l]_<<<{p_1} &
\ \ (\mathbb{P}^1)^{n-2}  \\
}
$$
o\`u $p_1$ et $p_2$ sont les deux projections. L'image r\'eciproque
d'un point de
$(\mathbb{P}^1)^{n-2}
$ par
$
\ph_{R_1+R_n}$ est isomorphe \`a $ (\mathbb{P}^{1})^{2}$, les
restrictions de $ \ph_{R_1}$ et $\ph_{R_n} $
\'etant les projections sur chacun des facteurs. En d'autres termes, une
fibre de $\ph_{R_n}  $ est envoy\'ee isomorphiquement par $  \ph_{R_1}$
sur le facteur
$\PP^1$ de $W_1$. Cela signifie que le degr\'e du   morphisme produit
$$\varphi_{n-1}\times  \ph_{R_n} :X\to \PP^1\times (\PP^1)^{n-1},$$
qui est le degr\'e de la restriction de $ \varphi_{n-1}$ \`a une fibre
g\'en\'erale de
$\ph_{R_n}$, vaut $1$. Ce morphisme \'etant fini, c'est un
isomorphisme et la preuve de la proposition
est achev\'ee.\finpreuve

\medskip

\section{La dimension quatre}
\marque Montrons l'in\'egalit\'e ($\ast$) en dimension $4$.
Une remarque essentielle est que sur une vari\'et\'e
de Fano de dimension $4$ et de pseudo-indice
$\ge2$, il n'y a pas
de petite contraction (\voir 1.2).

\medskip

\noindent {\bf Proposition \num} {\em
Soit $X$ une vari\'et\'e de Fano de dimension $4$.
Si $\iota_X =3$, on a $\rho_X \le 2$.
}

\medskip

\noindent {\em D\'emonstration.} Supposons $\rho_X \ge 3$.
Les fibrations extr\'emales sont de dimension relative au
moins $2$
(\voir 1.2) et il
y en a donc au plus $2$ (\cite{Wis91}, Theorem (2.2)).
Il existe ainsi au moins une contraction
divisorielle contractant son lieu exceptionnel sur un point
(\voir 1.2). Il en r\'esulte
l'existence sur $X$ d'une contraction lisse de centre de codimension
$2$ (\cite{Wis91}, Corollary (1.3)), ce qui est absurde puisqu'une
telle contraction est
de longueur $1$. \finpreuve

\medskip

\hs Le r\'esultat principal de cette section est le

\medskip

\noindent {\bf Th\'eor\`eme \num}
{\em Soit $X$ une vari\'et\'e de Fano de dimension $4$.
Si $\iota_X =2$, on a $\rho_X \le 4$.
}

\subsection{Quelques lemmes interm\'ediaires}
La d\'emonstration du th\'eor\`eme~3.3 repose sur les lemmes suivants.

\medskip

\noindent {\bf Lemme \num} {\em
Une vari\'et\'e de Fano $X$ de dimension $4$,
de pseudo-indice $2$ et de nombre de Picard
$\ge5$ n'a pas
de fibration extr\'emale $\varphi : X \to Y$ dont toutes
les fibres sont de dimension $1$.
}

\medskip

\noindent {\em D\'emonstration.} Raisonnons par l'absurde.
Le morphisme $\varphi : X \to Y$ est une fibration lisse car  $\iota_X
=2$, de sorte que $Y$ est une vari\'et\'e de Fano lisse de dimension $3$
d'apr\`es~1.3. Puisque $\rho_Y \ge 4$, il existe une
contraction lisse
$\pi : Y \to Z$ de centre une courbe lisse (\cite{MM81}, Theorem~5).
Toute courbe de $Y$ contract\'ee par
 $\pi$ est une courbe rationnelle de fibr\'e normal
dans $Y$ \'egal \`a ${\mathcal O} (-1)\oplus {\mathcal O} $, donc de
degr\'e anticanonique $1$. Ceci, avec le
lemme 2.5, contredit l'hypoth\`ese
$\iota_X =2$.\finpreuve

\medskip

\noindent {\bf Lemme \num} {\em
Soit $X$ une vari\'et\'e de Fano de dimension $4$,
de pseudo-indice $2$ et de nombre de Picard $\ge5$.
Toute contraction extr\'emale $ X \to Y$ est ou bien une
contraction divisorielle lisse de centre une courbe, autrement dit
est l'\'eclatement d'une vari\'et\'e lisse $Y$ le long
d'une courbe lisse de $Y$,
ou bien une fibration de dimension relative $1$.
}

\medskip

\hs Dans le cas d'une fibration, le
lemme~3.5 montre qu'il y a au moins une fibre
de dimension $2$.

\medskip

\noindent {\em D\'emonstration.}
Remarquons tout d'abord  que les fibrations extr\'emales
$X \to Y$ v\'erifient toutes $\dim(Y)=2$ ou $\dim(Y)=3$.
En effet, si $\dim(Y) = 0$ (resp. $\dim(Y) = 1$),
on a $\rho_X =1$ (resp. $\rho_X =2$).
D'autre part, comme $\dim(X) =4$ et $\rho_X \ge 5$,
il y a au moins une contraction extr\'emale
divisorielle (\cite{Wis91}, Theorem (2.2)).

\hs Soient $\varphi : X \to Y$ une telle contraction
et $E$ son diviseur exceptionnel.
Les fibres non triviales de $\varphi$ sont
de dimension $\ge2$ (\voir 1.2).
Si $\dim (\varphi (E)) =0$, il existe sur $X$ une
fibration extr\'emale dont toutes les fibres sont de dimension $1$ ou une
contraction extr\'emale birationnelle lisse de centre de codimension $2$
(\cite{Wis91}, Corollary (1.3)), ces deux situations \'etant exclues
respectivement par le lemme 3.5 et par l'hypoth\`ese $\iota_{X}=2$.
Ainsi, $\dim(\varphi (E)) =1$ et $\varphi$ s'identifie \`a l'\'eclatement
d'une vari\'et\'e lisse $Y$ le long de la courbe lisse $\varphi (E)$
(\cite{AW98}, Theorem 4.1).

\hs V\'erifions enfin qu'il n'existe pas de fibration
extr\'emale $\pi : X \to S$ o\`u
$S$ est une surface. Supposons qu'une telle fibration existe; $S$ est
alors une surface lisse (\cite{ABW92}, Proposition
1.4.1) et une fibre
 non triviale $F$ de $\varphi$ est  isomorphe
\`a
$\PP^2$. Puisque
$\varphi_{|F} : F \to S$ est fini, $S$ est \'egalement  isomorphe \`a
$\PP^2$ (\cite{Laz83}), de sorte que $\rho_X =2$, ce qui est
absurde.\finpreuve

\medskip

\hs Le r\'esultat suivant pr\'ecise le cas des \'eclatements de
centre une courbe lisse.

\medskip

\noindent {\bf Proposition \num} {\em
Soit $X$ une vari\'et\'e de Fano de dimension $n \ge 4$ et
de pseudo-indice $\iota_X \ge 2$.
Si $\pi : X \to Y$ est l'\'eclatement d'une
vari\'et\'e lisse $Y$ le long d'une courbe lisse, $Y$ est une vari\'et\'e
de Fano et $\iota_Y \ge \iota_X$. }

\medskip

{\em D\'emonstration.}
Si $Y$ n'est pas une vari\'et\'e de Fano,
le centre de l'\'eclatement $\pi$ est une courbe rationnelle lisse de
fibr\'e normal
${\mathcal O}_{\PP^1}(-1)^{\oplus n-1}$ (\cite{Wis91}, Proposition (3.5)),
de sorte que $\iota_X =1$, ce qui contredit l'hypoth\`ese.

\hs V\'erifions
ensuite l'in\'egalit\'e
$\iota_Y
\ge 2$. Soient $C$ le centre de
$\pi$ et $E$ son diviseur exceptionnel.
Si $C'$ est une courbe rationnelle de $Y$ distincte de $C$, et
si on note encore
$C'$ sa transform\'ee stricte dans $X$,
on a $-K_Y\cdot C'=-K_X\cdot C'+(n-2)E\cdot C' \ge \iota_X$.
Si $\iota_Y < \iota_X$, on ne peut donc avoir $-K_Y\cdot C'=\iota_Y$, de
sorte que $C$ est  une courbe
rationnelle qui satisfait $-K_Y\cdot C=\iota_Y$.

\hs Soit $N_{C/Y}=
\bigoplus_{i=1}^{n-1} {\mathcal O}_{\PP^1}(a_i)$ son fibr\'e normal. On a
d'une part
$$\iota_Y=-K_Y\cdot C =\sum_{i=1}^{n-1}a_i +2,$$ et d'autre part, si
$C_i\subset E$ est la courbe d\'efinie par le fibr\'e quotient
${\mathcal O}_{\PP^1}(a_i)$ de $N_{C/Y}$,
$$-K_X\cdot C_i = \iota_Y -(n-2)a_i \ge \iota_X.$$
On en d\'eduit $\iota_Y \ge \iota_X$, sauf si tous les $a_i$ sont
strictement n\'egatifs, ce qu'exclut l'\'egalit\'e $\sum_{i=1}^{n-1}a_i =
\iota_Y-2
\ge -1$.\finpreuve

\subsection{Vari\'et\'es de Fano  sp\'eciales}
Nous nous int\'eressons ici \`a certaines vari\'et\'es de Fano,
que nous appelons \og sp\'eciales\fg, faute d'une meilleure terminologie.

\medskip

\noindent {\bf D\'efinition.}
Une vari\'et\'e de Fano de dimension $n$ sera dite {\em sp\'eciale} si
toutes ses contractions extr\'emales sont ou bien
des \'eclatements lisses  de centre une courbe
lisse ou bien des fibrations
de dimension relative $1$ et si elle poss\`ede au moins
une contraction extr\'emale birationnelle.

\medskip

\noindent {\bf Exemple.} La vari\'et\'e $X = \PP^n_a \times \PP^1$,
o\`u $\PP^n_a$ d\'esigne l'\'eclatement de $\PP^n$
de centre $a$, est une vari\'et\'e de Fano sp\'eciale, de pseudo-indice
$2$ et de nombre de Picard $3$. En effet,
$X$ poss\`ede $3$ contractions extr\'emales:
$X \to \PP^n_a$ et $X \to \PP^{n-1}\times \PP^1$, fibrations en $\PP^1$,
et
$X \to \PP^n \times \PP^1$, contraction lisse
de centre $\{a\}\times \PP^1$.

\medskip

\hs Cet exemple est d'une certaine fa\c con extr\'emal comme
le montre le

\medskip

\noindent {\bf Th\'eor\`eme \num} {\em
Une vari\'et\'e
de Fano sp\'eciale $X$ de dimension $\ge 4$ v\'erifie
$\rho_X \le 3$.
}

\medskip

\hs La d\'emonstration de ce r\'esultat occupe la fin de ce
paragraphe.

\medskip

\noindent {\bf Lemme \num} {\em
Soit $X$ une vari\'et\'e de Fano sp\'eciale
de dimension $\ge 4$ et de nombre de Picard  $\ge4$.
\begin{enumerate}
\item[(a)]
La vari\'et\'e $X$ poss\`ede au plus une
fibration extr\'emale
de dimension relative $1$.
\item[(b)] Si $X$ poss\`ede une fibration
extr\'emale de dimension relative $1$, d'ar\^ete $R$, les contractions
extr\'emales birationnelles ont pour centre une courbe rationnelle \`a
fibr\'e normal trivial, leur diviseur exceptionnel $E$ v\'erifie $E \cdot
R =0$  et $\iota_X
\le 2$.
\end{enumerate}
}

\medskip

{\em D\'emonstration.}
Soient $\pi : X \to W$ une contraction extr\'emale birationnelle
de centre une courbe
lisse $C \subset W$ et de diviseur exceptionnel
$E$ et
$\varphi : X \to Y$ une fibration extr\'emale d'ar\^ete
$R$ et de dimension relative $1$.
Nous allons montrer que $\pi$ d\'etermine  $R$, donc $\ph$.

\hs Le morphisme $\varphi_{|E} : E \to Y$
n'est pas surjectif car $\rho_E =2$ et $\rho_Y \ge 3$;
comme l'intersection d'une fibre de $\varphi_{|E}$ avec une fibre
de $\pi_{|E}$ est finie, chaque fibre de $\varphi_{|E}$ est de dimension
$1$. En particulier, $\ph(E)$ est de codimension $1$ dans $Y$;
l'image r\'eciproque par   $\ph$ d'un point
g\'en\'eral de
$\ph(E)$ est  de dimension $1$, donc co\"incide avec
 son image r\'eciproque par
$\varphi_{|E}$. C'est en particulier une courbe rationnelle qui
domine $C$, de sorte que cette derni\`ere aussi est rationnelle.

\hs
Le diviseur $E$ est isomorphe \`a $\PP (N_{C/W}^{\ast})$,
o\`u $N_{C/W}=\bigoplus_{i=1}^{n-1} {\mathcal O}_{\PP^1}(a_i)$.
La fibration induite $\widetilde\varphi_{|E}  : E \to
\widetilde\varphi(E)$, o\`u $\widetilde\varphi(E)$ est la
normalisation de
$\varphi(E)$, est
\'equidimensionnelle, de fibre
g\'en\'erale $\PP^1$; c'est  une contraction extr\'emale
puisque $\rho_{E}=2$. Le c\^one de Mori de $\PP (N_{C/W}^{\ast})$
est engendr\'e par
deux courbes: une droite d'une fibre de $\pi_{|E} : E\to C$
et une section de $\pi_{|E} $. Il s'ensuit
que $E$ est isomorphe \`a $ \PP^{n-2} \times \PP^1$, que
$\pi_{|E}$ et $\varphi_{|E}$ sont les deux projections et
 que $R$ est engendr\'e par une courbe du type
$\{*\}\times \PP^1 $. Comme $n\ge4$,
il y a donc au plus une fibration extr\'emale
de dimension relative $1$.

\hs Ce qui pr\'ec\`ede montre aussi que
si $X$ poss\`ede une fibration
extr\'emale de dimension relative $1$, les contractions
extr\'emales birationnelles ont pour centre une courbe  rationnelle de
fibr\'e normal trivial. En effet, on a vu que $\PP (N_{C/W}^{\ast})$ est
isomorphe
\`a
$\PP^{n-2} \times \PP^1$, de sorte que
$N_{C/W}$ est isomorphe \`a $ {\mathcal O}_{\PP^1}(a)^{\oplus n-1}$
pour un entier $a$ convenable et, puisque les courbes
du type $\{*\}\times \PP^1$ ont un fibr\'e normal (dans
$X$) trivial (ce sont les fibres d'une fibration), c'est que $a=0$.
\finpreuve

\medskip

\noindent {\bf Lemme \num} {\em
Soit $X$ une vari\'et\'e projective lisse de dimension $4$
et soit $\pi : X \to W$ (resp. $\pi' : X \to W'$, $\pi'' : X \to W''$)
une contraction extr\'emale birationnelle lisse
de centre une courbe
lisse $C \subset W$ (resp. $C' \subset W'$, $C'' \subset W''$)
d'ar\^ete $R$ (resp. $R'$, $R''$) et de diviseur exceptionnel
$E$ (resp. $E'$, $E''$).
\begin{enumerate}
\item[(a)] Si $R \neq R'$ et $E\cap E' \neq \emptyset$, on a $E \cdot R'
>0$ et
$E' \cdot R >0$.
\item[(b)] Si $R$, $R'$ et $R''$
sont deux \`a deux distinctes et si
$E\cap E' \neq \emptyset$ et $E\cap E'' \neq \emptyset$, on a $E'\cap
E''\neq
\emptyset$.
\end{enumerate}
}

\medskip

{\em D\'emonstration.}
V\'erifions le point (a). Soit $x$ un point de   $E\cap E'$. La surface
$\pi^{-1}(\pi (x)) $, isomorphe \`a $ \PP^2$, n'est pas contenue
dans $E'$ car $R \neq R'$. Il existe donc
une courbe dans $\pi^{-1}(\pi (x))$ passant par
$p$ et non contenue dans $E'$,  d'o\`u (a) en \'echangeant
les r\^oles de $E$ et $E'$.

\hs V\'erifions le point (b). Soit $y \in C$. Les courbes
$\pi^{-1}(y) \cap E'$ et $\pi^{-1}(y) \cap E''$
sont contenues dans $\pi^{-1}(y) $, qui est isomorphe \`a $ \PP^2$, de
sorte que
$\pi^{-1}(y) \cap E' \cap E''\ne \emptyset$,  d'o\`u (b).\finpreuve

\medskip

\subsection{D\'emonstration du th\'eor\`eme 3.9}
Supposons $\rho_X \ge 4$. Par hypoth\`ese, il existe une contraction
extr\'emale birationnelle $ \pi: X \to W$ d'ar\^ete $R_1$ et de diviseur
exceptionnel
$E_1$, d'image une vari\'et\'e de Fano lisse (proposition 3.7) et de
centre une courbe lisse
$C \subset W$.

\hs Soit $C_1$ une courbe telle
que $E_1 \cdot C_1 >0$.
Cette courbe est nu\-m\'e\-ri\-quement combinaison
lin\'eaire \`a coefficients
rationnels strictements positifs de classes engendrant des ar\^etes.
L'une d'entre elles, not\'ee $R_2$, satisfait $E_1 \cdot R_2 > 0$.
Le lemme 3.10(b) entra\^ine que la contraction $\ph_{R_2}$ est
birationnelle.

\hs Notons $R_1,R_2,\dots,R_k,R_{k+1},\dots,R_m,R_F$
les ar\^etes de $\NE (X)$, o\`u
\begin{itemize}
\item[$\bullet$] les contractions $\ph_{R_i}$ sont birationnelles de
diviseur exceptionnel $E_i$ (noter que $\pi=\ph_{R_1}$);
\item[$\bullet$]
$\ph_{R_F}$ est l'\'eventuelle unique fibration et,  d'apr\`es le lemme
3.10(b),
$E_i \cdot R_F =0$ pour tout $1\le i \le m$;
\item[$\bullet$]
$E_1  \cap  E_i \neq \emptyset$ pour
$2\le i \le k$ d'o\`u, d'apr\`es le lemme 3.11(a),
$E_1 \cdot R_i >0$;
\item[$\bullet$]
 $E_1 \cap E_i = \emptyset$ pour $k+1\le i \le m$, d'o\`u, d'apr\`es le
lemme 3.11(b),
$E_2 \cap E_i = \emptyset$.
\end{itemize}

\hs Supposons $n\ge5$.  L'intersection
$E_1 \cap E_2$ n'est pas vide, de
dimension au moins $n-2 \ge 3$: il y a donc des courbes de
$E_1 \cap E_2$ contract\'ees \`a la fois
par $\ph_{R_1}$ et par $\ph_{R_2}$, ce qui est absurde.

\hs Supposons $n=4$ et v\'erifions que $C$ est extr\'emale dans $\NE(W)$.
Rappelons que
ce c\^one est engendr\'e par
$\pi(R_F)$ et les
$\pi(R_i)$ pour $i\ge 2$.
Pour $2\le i \le k$,   la classe de $C$
appartient \`a $\pi(R_i)$,  ainsi qu'\`a
$\pi(R_F)$: en effet, chaque fibre
non triviale de $\ph_{R_i}$ intersecte $E_1$ le long d'une courbe
$C_i$ telle que $\pi(C_i)=C$, de m\^eme, toute fibre $F$ de
$\ph_{R_F}$ intersectant $E_1$ est contenue dans $E_1$
et satisfait donc $\pi(F)=C$.
Si la classe de $C$ n'est pas extr\'emale,
$\NE (W)$ est donc engendr\'e par les ar\^etes $\pi(R_i)$
pour $i\ge k+1$.
L'image $\pi(E_2)$ est un diviseur effectif de $W$
num\'eriquement trivial puisque $E_2 \cdot R_{F}=E_2 \cdot R_i =0$,
ce qui est absurde.

\hs \'Etudions la contraction extr\'emale $\varphi : W \to Y$
d'ar\^ete engendr\'ee par $C$ et montrons
que le diviseur $\pi(E_2)$ est contract\'e sur un point par
$\varphi$. Notons $W_2$ l'image de  $\ph_{R_2}$ et $C_2\subset W_2$
son centre.
Pour tout  $y\in C_2$, la fibre $\ph_{R_2}^{-1}(y)$, isomorphe \`a
$\PP^2$, contient une courbe dont l'image par $\pi$ est $C$ et qui
est donc contract\'ee par $\varphi$ sur le point
$\varphi(C)$. D'autre part, $\pi(\ph_{R_2}^{-1}(y))$ est
\'egalement  contract\'e sur le point $\varphi(C)$.
La contraction
$\varphi : W \to Y$ est donc ou bien divisorielle ou bien une
fibration. Soit $C' \in R_2$ une courbe rationnelle de $X$ telle que
$E_2\cdot C' =-1$. Il existe un entier $r\ge 1$ tel que
$$ \pi(E_2)\cdot \pi_*(C') = \pi^*(\pi(E_2))\cdot C'
=  (E_2 + rE) \cdot C' = -1+r E\cdot C'.$$
Or $E\cdot C' \ge 1$; on en d\'eduit   $\pi(E_2)\cdot \pi_*(C')\ge 0$
puis, puisque les courbes $C$ et $\pi_*(C')$ sont num\'eriquement
proportionnelles,
$\pi(E_2)\cdot C \ge 0$.\\
\hs Comme $\pi(E_2)$ est contract\'e sur un point
par $\varphi $, le calcul d'intersection pr\'ec\'edent montre
que $\varphi $ est une fibration. En particulier, $\dim(Y) \le 1$
et
$\rho_W \le 2$, ce qui est absurde puisqu'on a suppos\'e $\rho_X \ge 4$.
\finpreuve

\subsection{D\'emonstration du th\'eor\`eme~3.3}
Le lemme~3.6 montre que si $X$ est une vari\'et\'e de Fano de dimension
$4$, de pseudo-indice $2$ et de nombre de Picard
 $\ge5$,   la contraction extr\'emale
associ\'ee
\`a toute ar\^ete extr\'emale de $X$ est soit une contraction
divisorielle lisse de centre une courbe lisse,
soit une fibration de dimension relative $1$.
De plus, comme $\rho_X \ge 5$, il existe sur $X$
au moins une contraction extr\'emale birationnelle
(\cite{Wis91}, Theorem (2.2)).
C'est donc que $X$ est sp\'eciale et le th\'eor\`eme~3.9 permet de
conclure \`a une absurdit\'e.\finpreuve

\subsection{Familles propres de courbes rationnelles.}
Nous renvoyons au livre \cite{Kol99}
pour plus de d\'etails sur les
notations et les rappels qui suivent. Soit $X$
une vari\'et\'e complexe, projective, lisse et connexe.
Soit $\textup{Hom}_{\textup{bir}}(\PP^1,X)$ le sch\'ema des morphismes
birationnels de $\PP^1$ vers $X$ et soit
$\textup{Hom}_{\textup{bir}}^{n}(\PP^1,X)$
sa normalisation. Le groupe lin\'eaire $\textup{PGL}(2,\CC)$ agit
sur $\textup{Hom}_{\textup{bir}}^{n}(\PP^1,X)$
et $\textup{Hom}_{\textup{bir}}^{n}(\PP^1,X)\times\PP^1$.
Les quotients g\'eom\'etriques au sens de Mumford existent
et seront respectivement not\'es
$\textup{RatCurves}^{n}(X)$ et $\textup{Univ}^{rc}(X)$.
Soit $V\subset\textup{RatCurves}^{n}(X)$ une famille propre irr\'eductible
de courbes rationnelles irr\'eductibles sur $X$ et soit
$\mathcal U\subset\textup{Univ}^{rc}(X)$ la famille universelle
\begin{equation*}
\begin{CD}
\mathcal U @){\ev})) X \\
@V{\pi}VV \\
V
\end{CD}
\end{equation*}
Notons $\lieu (V) = \ev ({\mathcal U})$ l'ensemble
des points de $X$ par lesquels il passe une courbe
rationnelle $C$ de $X$
telle que $[C] \in V$. Soit $x\in\lieu (V) $; on note
$V_x=\pi ( \ev^{-1}(x) ) \subset V$
les courbes de $V$ passant par $x$, puis ${\mathcal U}_x = \pi^{-1}(V_x)$
et $\lieu (V_x) = \ev ({\mathcal U}_x)$. Les dimensions de ces
diff\'erentes  vari\'et\'es satisfont
\begin{eqnarray}\label{dim}
\dim (V) &\ge& -K_X\cdot V + n -3 \\
 \dim (V_x)&\ge& -K_X\cdot V -2\nonumber
\end{eqnarray}
o\`u $-K_X\cdot V$ d\'esigne l'intersection
$-K_X \cdot C$ pour une courbe rationnelle $C$ de $X$
telle que $[C]\in V$. Si $x$ est un point g\'en\'eral de $\lieu (V)$, on a
\begin{equation}\label{lieu}
 \dim (\lieu (V)) + \dim (V_x)= \dim(V) +1
\end{equation}
et, par le lemme de cassage (\cite{cass}, Theorem 6),
$$\dim (\lieu (V_x) ) = \dim (V_x) +1 \ge -K_X\cdot V -1.
$$
 Nous utiliserons de fa\c con r\'ep\'et\'ee le r\'esultat suivant
(\cite{Kol99}, II.4.21).

\medskip

\noindent {\bf Lemme \num} {\em
Soit $V \subset\textup{RatCurves}^{n}(X)$ une famille propre et
irr\'eductible de courbes rationnelles irr\'eductibles sur une vari\'et\'e
$X$ projective et lisse et soit $x\in \lieu (V)$. Toute
courbe trac\'ee sur $\lieu (V_x) $ est num\'eriquement
proportionnelle \`a une courbe $C$ telle que $[C] \in V$.
}

\medskip

\marque La famille $V$ d\'etermine une relation d'\'equivalence sur $X$
pour  laquelle des points $x$ et $x'$ de $X$ sont \'equivalents s'il
existe une cha\^ine connexe de courbes rationnelles de $V$ passant par
$x$ et $x'$. Il existe un ouvert
$X_{0}\subset X$ et un morphisme propre $X_{0}\to Z_{0}$
\`a fibres connexes vers une vari\'et\'e normale
dont les fibres
sont des classes d'\'equivalence pour la relation
pr\'ec\'edente (\cite{Kol99}, IV.4.16).

\subsection{Les cas d'\'egalit\'e.} L'objet de ce paragraphe
est l'\'etude des cas
d'\'egalit\'e dans $(\ast)$.

\medskip

\noindent {\bf Proposition \num} {\em Soit $X$ une vari\'et\'e de Fano de
dimension $4$. Si $\rho_X(\iota_{X}-1)=4$, la vari\'et\'e $X$ est
isomorphe \`a $(\mathbb{P}^{\iota_{X}-1})^{\rho_{X}}.$
}

\medskip

\noindent {\em D\'emonstration.}
Si $\rho_{X}=1$ et $\iota_{X}=5$, la vari\'et\'e $X$ est isomorphe \`a
 $\PP^4$ (\cite{CMS00}, \cite{Keb01}).

\smallskip

\hs Supposons $\rho_{X}=2$ et $\iota_{X}=3$.
Notons $\ph_{1}:X\to W_{1}$
et $\ph_{2}:X\to W_{2}$ les deux contractions extr\'emales,
d'ar\^etes respectives $R_{1}$ et $R_{2}$.
Elles sont ou bien divisorielles ou bien des fibrations (\voir 1.2).
Supposons par exemple que $\ph_{1}$ soit divisorielle.
Le lieu exceptionnel $E_{1}$ de $\ph_{1}$ est contract\'e sur un
point  par $\ph_{1}$ (\voir 1.2). Rappelons que l'intersection de
fibres de deux contractions extr\'emales diff\'erentes est finie.
Les fibres non triviales de $\ph_{2}$ sont de dimension au moins $2$, donc
ne rencontrent pas $E_{1}$, autrement dit, $\ph_{2}$ est
\'egalement  divisorielle de
lieu exceptionnel $E_{2}$ disjoint de $E_{1}$. Le diviseur $-E_{1}$ est donc
num\'eriquement effectif puisque
$-E_{1}\cdot R_{1}>0$ et $-E_{1}\cdot R_{2}=0$,
ce qui est manifestement absurde.
Les morphismes $\ph_{1}$ et $\ph_{2}$ sont
donc des fibrations \'equidimensionnelles de dimension relative $2$
(\voir 1.2) et
$W_{1}$ et $W_{2}$ sont lisses (\cite{ABW92}, Proposition 1.4.1).
Une fibre g\'en\'erale $F_{1}$ (resp. $F_{2}$)
de $\ph_{1}$ (resp. $\ph_{2}$) est de pseudo-indice
$3$
par la
formule d'adjonction
et donc isomorphe \`a $\PP^{2}$.
Le morphisme $F_{1}\to W_{2}$
(resp. $F_{2}\to W_{1}$) est fini et $W_{1}$
(resp. $W_{2}$)
est isomorphe
\`a $\PP^{2}$ par le th\'eor\`eme de Lazarsfeld (\cite{Laz83}).\\
\hs Les fibres de $\ph_{2}$ ne sont {\em a priori}
pas des sections de $\ph_{1}$. Nous allons
montrer que de telles sections existent.
Si $\ell\subset\mathbb{P}^{2}$ est une droite g\'en\'erale,
$X_{\ell}=\ph_{1}^{-1}(\ell)$
est lisse et connexe. Les fibres
g\'en\'erales de $X_{\ell}\to \ell$ sont isomorphes
\`a $\mathbb{P}^{2}$.
Il existe donc une section de $X_{\ell}\to \ell$ par le
th\'eor\`eme de Tsen (\cite{Kol99}, Theorem IV.6.5).

\hs Soit donc
$C\subset X$ une courbe rationnelle v\'erifiant
$C\cdot \ph_{1}^{*}\mathcal{O}_{\PP^{2}}(1)=1$, de degr\'e minimal
relativement
au diviseur ample $-K_{X}$. La courbe $C$ d\'etermine une famille
{\em propre} irr\'eductible $V\subset\textup{RatCurves}^{n}(X)$ de courbes
rationnelles irr\'eductibles sur $X$ (\voir 3.14).
La dimension de $V$ est au moins $-K_X\cdot C + 1\ge 4$ et celle
de $V_{x}$, pour $x\in \lieu (V)$,
au moins $-K_X\cdot C - 2\ge 1$.
Si la dimension de $\lieu (V_{x})$ est au moins $3$, il rencontre une
fibre g\'en\'erale de $\ph_{1}$ au moins le  long d'une courbe et $[C]\in
R_{1}$ (\voir lemme 3.15) ce qui est absurde  par le choix de $C$. Ainsi,
$\dim( \lieu (V_{x}))=2$ et le lieu de $V$ est $X$ par
la formule (\ref{lieu}) de 3.14.

\hs Il existe un ouvert
$X_{0}\subset X$, une vari\'et\'e normale $Z_0$ et un morphisme propre
$q:X_{0}\to Z_{0}$
\`a fibres connexes dont les fibres  sont des classes
  pour la relation  d'\'equivalence
d\'etermin\'ee par $V$ (\voir 3.14). La famille $V$ est couvrante et
$Z_{0}$
est donc
de dimension au plus $3$. Soient $F$ une fibre g\'en\'erale de
$q$ et
$x$ un point de $F$.  Comme $\dim( \lieu (V_{x}))=2$, on a $\dim (F) \ge
2$. Si $F$ est de dimension au moins $3$,
elle rencontre une fibre
g\'en\'erale de $\ph_{1}$ au moins le  long d'une courbe.
Comme $F$ est couverte par une famille propre de courbes
rationnelles irr\'eductibles de $V$, toutes les courbes de $F$ sont
alg\'ebriquement \'equivalentes \`a un multiple de $C$ (\cite{Kol99}, IV
3.13.3)
 et $[C]\in
R_{1}$,
ce qui est
\`a nouveau absurde. Ainsi $F$ est une surface de pseudo-indice $3$ par
la formule d'adjonction, donc isomorphe \`a $\mathbb{P}^{2}$; on a
$F=\lieu (V_{x})$ pour tout
$x$ dans $F$, les courbes de $V_x$ \'etant les droites passant par $x$.
Notons que
$F$ est une section de
$\ph_1$.\\
\hs V\'erifions que l'on peut prendre $X_{0}=X$.
Le fibr\'e normal \`a $F $ dans $X$ est trivial et le sch\'ema de
Hilbert de $X$ est donc lisse au point $[F]$. Soient $H$ son unique
composante passant par $[F]$ et $X'\subset X\times H$ la famille
universelle. Quitte \`a remplacer $Z_{0}$ par un ouvert dense, on
peut supposer que $q$ est plat et on a un diagramme
\begin{equation*}
\begin{CD}
X_0 @)))X' @){\varphi})) X \\
@V{q}VV@V{\pi}VV \\
Z_0 @)))H
\end{CD}
\end{equation*}
o\`u le carr\'e est cart\'esien et o\`u la compos\'ee des fl\`eches
horizontales sup\'erieures est l'inclusion. Le morphisme
$\pi$ s'identifie, au-dessus d'un ouvert
$H_{0}$ convenable de $H$, au morphisme $X_{0}\to Z_{0}$.

\hs Fixons un point
$t_{0}$ de $ H$ et des points $x'_{1}$ et $x'_{2}$ de $\pi^{-1}(t_0)$.
Notons
$T\to H$ un germe de courbe lisse passant par $t_{0}$ et
rencontrant $H_0$. Le sch\'ema  $\pi_{T}:X'_{T}\to T$
obtenu
par le changement de base $T\to H$ est irr\'eductible: il
existe donc des courbes $T_1\to X'_T$ et $T_2\to X'_T$ dominant
$T$ et rencontrant $\pi_T^{-1}(t_0)$ uniquement en
$x'_1$ et $x'_2$ respectivement. Pour tous points
$t_1\in T_1$ et $t_2\in T_2$ situ\'es au-dessus du m\^eme point
g\'en\'eral $t$ de
$T$, les points correspondant de $X$ peuvent \^etre joints par une
(unique) courbe $C_{t_1,t_2}\subset\pi_{T}^{-1}(t)$ de  $V$. Cette
famille \'etant propre, les courbes $C_{t_1,t_2}$  d\'eg\'en\`erent
vers une courbe de $V$, qui est  en particulier irr\'eductible, joignant
$x'_1$ et $x'_2$. Les fibres de
$\pi$ sont donc  irr\'eductibles et  deux points
quelconques sont reli\'es par une courbe de $V$.

\hs Supposons que le morphisme $\varphi$ contracte
 une courbe irr\'eductible
$C\subset X'$  vers un point $x$ de $X$.
Celle-ci est horizontale
pour $\pi$ et le lieu couvert par les $\varphi(X'_{t})$, pour
$t\in\pi(C)$, est de dimension $3$.
Ainsi, le lieu couvert par les courbes de $V_{x}$ est de dimension $3$,
ce qui est exclu par
les arguments pr\'ec\'edents. Le morphisme $\ph$ est donc
birationnel
et fini:
c'est un isomorphisme.

\hs Notons $H'\to H_{\textup{r\'ed}}$ la normalisation; il
existe une factorisation $\pi:X'\stackrel{\pi'}{\to}H'\to
H_{\textup{r\'ed}}$ et le morphisme $q'=\pi'\circ\varphi^{-1}:X\to H'$
\'etend  $q$. Le morphisme produit $\ph_{1}\times q':X\to
\PP^{2}\times H'$ est birationnel
(puisqu'une fibre g\'en\'erale de $q'$ est une section de $\ph_1$);
puisque $\rho_X=2$, il est fini: c'est un isomorphisme. On en d\'eduit
que $H'$ est lisse et, comme les fibres  de $\ph_1$, isomorphe
\`a $\PP^{2}$.

\smallskip

\hs Supposons $\rho_{X}=4$ et $\iota_{X}=2$. Par 1.2, toute
contraction extr\'emale birationnelle est divisorielle et ses fibres non
triviales sont de dimension $\ge2$. V\'erifions que les fibrations
extr\'emales,  s'il en existe, sont de dimension relative $1$. Soit
$\varphi:X\to Z$ une telle contraction; supposons qu'elle soit de
dimension relative au moins $2$.  Notons que $Z$ est lisse (\cite{ABW92},
Proposition 1.4.1). Puisque $\rho_{X}=4$, il existe (\cite{Wis91},
Theorem (2.2)) une contraction extr\'emale birationnelle; ses fibres non
triviales sont de dimension $\ge2$, de sorte que
$\varphi$ est
de dimension relative $2$. La vari\'et\'e $W$ est donc lisse
et $\psi$ s'identifie
\`a l'\'eclatement d'une courbe lisse dans $W$ (\cite{AW98}, Theorem
4.1). Toute fibre de $\psi_{|E}:E\to\psi(E)$ est en
particulier
isomorphe \`a
$\mathbb{P}^{2}$ et domine $Z$. Il en r\'esulte que $Z$ est isomorphe
\`a
$\mathbb{P}^{2}$ (\cite{Laz83})
ce qui est absurde puisque $\rho_Z=\rho_{X}-1=3$.\\
\hs La vari\'et\'e $X$ n'est pas sp\'eciale (th\'eor\` eme 3.9)
et ou bien toutes les contractions \'el\'e\-men\-taires
sont des fibrations de dimension relative $1$, ou bien il existe une
contraction birationnelle $\psi:X\to W$, divisorielle, dont le lieu
exceptionnel $E$ est contract\'e sur un point. V\'erifions
que ce dernier cas ne peut pas se produire.
Il existe alors une contraction extr\'emale $\varphi:X\to Z$
dont le lieu exceptionnel rencontre
$E$. L'intersection d'une fibre de $\ph$ avec $E$ \'etant finie, cette
fibre est de dimension au plus $1$.
Puisque $\iota_{X}=2$, le morphisme $\ph$ est une fibration lisse
de
dimension relative $1$ par 1.3. On
en d\'eduit que $Z$ est une vari\'et\'e de Fano de dimension $ 3$, de
nombre de  Picard $ 3$ et de
pseudo-indice au moins $2$ (lemme 2.5) et donc exactement $2$ (\voir
section~2).
Finalement, $Z$ est isomorphe
\`a $(\mathbb{P}^{1})^{3}$ (proposition 2.3).
Posons $\ell =\PP^1\times\{t,t'\}$, pour $t$ et $t'$
g\'en\'eraux dans $\PP^1$. La surface
$X_{\ell}=\varphi^{-1}(\ell)$  est une surface de
Hirzebruch et la  courbe $\ell'=E\cap X_{\ell}$ est exceptionnelle car
contract\'ee par $\psi$: c'est donc une section de
$\varphi_{|X_{\ell}}:X_{\ell}\to\ell$.
Le morphisme $E\to Z$ est birationnel et fini: c'est
un isomorphisme.
On a $-K_{E}\cdot\ell'=2$ et
$-E\cdot\ell'>0$ puisque $\ell'$ est exceptionnelle dans $X_{\ell}$.
Or, par la
formule d'adjonction, on a
$-K_{E}\cdot\ell'=-K_{X}\cdot\ell'-E\cdot\ell'>2$
ce qui est absurde.\\
\hs En conclusion, toutes les contractions extr\'emales de $X$
sont des fibrations. Le r\'esultat
cherch\'e est donn\'e par la proposition 2.4.\finpreuve

\section{Le cas torique}

\hs Nous montrons dans ce paragraphe que si $X$ est
une vari\'et\'e de Fano \emph{torique} de dimension $n$, on a
$\rho_X (\iota_X -1) \le n$ lorsque
$\iota_X \ge \frac{1}{3}n+1$ ($n$ arbitraire) ou  $n\le 7$.

\subsection{Pr\'eliminaires.}
Soit $X$ une vari\'et\'e torique projective et lisse d'\'eventail
$\fx$. Nous renvoyons \`a \cite{Ful93} ou \cite{Oda88} pour les fondements de
la g\'eom\'etrie torique. Soit $G(\fx)$ l'ensemble des g\'en\'erateurs
primitifs
des c\^ones
de dimension $1$ dans $\fx$; son cardinal est $\rho_X+\dim (X)$.
Notons
$V(\sigma)$ l'adh\'erence de l'orbite correspondant \`a un \'el\'ement
$\sigma$ de $\fx$.

\hs Rappellons (\cite{bat1}, \S\ 2 et \cite{bat2}, \S\ 2.1)
qu'une \emph{collection primitive}
$P$ est un sous-ensemble de $G(\fx)$ minimal qui
n'engendre pas un c\^one de $\fx$. \`A toute collection primitive
$P=\{x_1,\dots,x_h\}$
est associ\'ee sa
\emph{relation primitive}
$$ x_1+\dots+x_h=a_1 y_1+\dots+a_k y_k, $$
o\`u $\langle y_1,\dots ,y_k\rangle$ est le plus petit c\^one dans $\fx$
contenant le point $x_1+\dots+x_h$ et o\`u les $a_i$ sont des entiers
strictement positifs.
Le degr\'e de $P$ est, par d\'efinition, $\deg (P)=h-\sum_ia_i$, et son
ordre est $|P|=h$.

\hs Le groupe $\Nun_{1}(X)$ des $1$-cycles
sur $X$ modulo \'equivalence num\'erique
s'identifie  au groupe des relations
entre les \'el\'ements de
$G(\fx)$; toute relation $\sum_{x\in G(\fx)} a_x x=0$
s'identifie \`a la classe d'\'equivalence num\'erique des
$1$-cycles
dont l'intersection
avec le diviseur $V(\langle x\rangle )$ est $a_x$.\\
\hs En particulier, la relation primitive
$x_1+\dots+x_h-(a_1 y_1+\dots+a_k y_k)=0$
d\'efinit un \'el\'ement $r(P)$ de $ \Nun_{1}(X)$. Cette classe est
toujours la classe d'un cycle effectif
et  $\deg (P)=-K_X\cdot
r(P)$  (\cite{bat1}, Theorem 2.15,
\cite{contr}, Lemma 1.4).
La vari\'et\'e $X$ est une
vari\'et\'e de Fano si et seulement si toute relation primitive est de
degr\'e strictement positif et  le pseudo-indice
$\iota_X$ v\'erifie alors
$|P|=h\ge\deg (P)\ge
\iota_X$.\\
\hs Une collection primitive $P$, ou sa relation associ\'ee $r(P)$, sont
dites  {\em contractibles} s'il existe une
application \'equivariante
$\ph : X\to W$ vers une vari\'et\'e torique compl\`ete $W$ telle que les
courbes irr\'eductibles contract\'ees par $\ph$ soient exactement
celles dont la classe est dans
$\mathbb{Q}_{\ge 0}r(P)$
(\cite{contr}, Definition 2.3). Avec les notations
pr\'ec\'edentes, le lieu exceptionnel de $\ph$ est $A=V(\langle
y_1,\dots,y_k\rangle)$ et
$\ph_{|A}: A\rightarrow \ph(A)$ est une fibration \'equivariante lisse
en $\PP ^{h-1}$. Si $\deg (P )< 2\iota_X$, la collection $P$ est
contractible (\cite{contr}, Theorem 4.1).

\subsection{Fibr\'es projectifs.} Soit $X$ une vari\'et\'e
torique projective et lisse. La contraction $\ph : X\to
Y$ d'une ar\^ete
num\'eriquement effective de $\NE(X)$
  est un
fibr\'e en espaces projectifs, $Y$ est une vari\'et\'e torique
et $\varphi$ est \'equivariant.

\hs Le lemme 2.5 permet de traiter facilement le cas des
vari\'et\'es toriques $X$ pour lesquelles toutes les collections
primitives de $\fx$
sont disjointes. Cette condition combinatoire est \'equivalente
\`a l'existence d'une suite $X=X_1\to X_2\to \dots\to X_r$ o\`u
$X_r$ est un espace projectif et $X_i$
est une fibration en espaces projectifs sur $X_{i+1}$ pour
$i=1,\dots,r-1$ (\cite{bat1}, Corollary 4.4). Rappelons enfin que
les vari\'et\'es toriques de nombre de Picard $2$
sont pr\'ecis\'ement les fibr\'es en espaces projectifs
sur l'espace projectif (\cite{Kle88}, Theorem 1)
et qu'elles satisfont donc aux
hypoth\`eses du

\medskip

\noindent {\bf Corollaire \num} {\em
Soit $X$ une vari\'et\'e de Fano torique. Si toutes les collections
primitives de
$\fx$ sont disjointes, on a $\rho_X(\iota_X-1)\le\dim (X)$; de plus, si
$\rho_X(\iota_X-1)=\dim (X)$, on a $X\simeq (\PP ^{\iota_X-1})^{\rho_X}$.
}

\medskip

{\em D\'emonstration.} Il existe donc des fibrations
$$ X=X_1\rightarrow X_2\rightarrow\dots\rightarrow X_r=\PP^{s_r}$$
o\`u $X_i$ est une fibration en espaces projectifs sur $X_{i+1}$
de fibre
$\PP ^{s_i}$ pour $i=1,\dots,r-1$. En appliquant le lemme
2.5,
on obtient,
par r\'ecurrence sur   $ i $, l'in\'egalit\'e
$\rho_X(\iota_X-1)\le\dim (X)$.
Supposons $\rho_X(\iota_X-1)=\dim (X)$. Comme $\rho_X=r$, $\dim
(X)=s_1+\dots+s_r$ et $\iota_X-1\le\min\{s_1,\dots,s_r\}$, on a
$$ s_1+\dots+s_r=r(\iota_X-1)\le r\min\{s_1,\dots,s_r\}\le
s_1+\dots+s_r, $$
et donc $s_1=\dots=s_r=\iota_X-1$.
On v\'erifie enfin que $X$ est un produit d'espaces projectifs
par r\'ecurrence
descendante sur l'entier $1\le i\le r-1$.
\finpreuve

\subsection{Le cas  $3\iota_X\ge \dim (X)+ 3$.}
Soit $X$ une vari\'et\'e de Fano torique. La preuve du r\'esultat
principal repose sur les deux lemmes suivants.

\medskip

\noindent {\bf Lemme \num} {\em
Soit $A\subset X$ une sous-vari\'et\'e irr\'eductible
invariante de codimension $k$. Si $k\le \iota_X-2$, alors $A$ est une
vari\'et\'e de Fano de nombre de Picard $\rho_X$ et  de pseudo-indice
$ \ge \iota_X-k$. }

\medskip

{\em D\'emonstration.}
Il suffit de prouver le lemme pour $k=1$
et $\iota_X\ge 3$: le r\'esultat g\'en\'eral
s'ensuit par r\'ecurrence sur $k$.
Soit donc $A=V(\langle x\rangle)$ un diviseur
irr\'eductible invariant et supposons $\iota_X\ge 3$.
Si $P$ est une collection
primitive dans $\Sigma_A$, ou bien $P$ est une collection
primitive dans $\fx$ et
$\deg_A (P)\ge \deg_X (P)\ge \iota_X$, ou bien
$P\cup\{x\}$ est une collection primitive dans $\fx$,
et $\deg_A (P)= \deg_X (P)-1\ge \iota_X-1$.
Comme $\iota_X\ge2$, le diviseur $A$ est une vari\'et\'e de Fano de
pseudo-indice
$\ge \iota_X-1$.
Puisque $\iota_X \ge 3$, il n'y a pas de collections primitives d'ordre $2$
(\voir 4.1): $\Sigma_A$ a donc exactement
un c\^one de dimension 1 de moins que $\fx$ et $\rho_{A}=\rho_{X}$.\finpreuve

\medskip

\noindent {\bf Lemme \num} {\em
Soit $X$ une vari\'et\'e de Fano torique de dimension $n$,
de pseudo-indice $\iota_X$ et de nombre de Picard $\rho_X$.
Si $4 \iota_X> n+4$, il existe dans $X$ une collection primitive
contractible telle que la codimension du lieu
exceptionnel de la contraction associ\'ee
soit $\le\iota_X-2$.
}

\medskip

{\em D\'emonstration.}
Supposons que toute relation primitive
contractible de $X$ soit de la forme
$ x_1+\dots+x_h=a_1 y_1+\dots +a_k y_k $
avec $k\ge \iota_X-1$. On a alors $h-k\ge h-\sum_i a_i \ge \iota_X$,
d'o\`u
$ h\ge \iota_X+k\ge 2\iota_X-1$. D'autre part, si $P$ est une
collection primitive non contractible,
on a $|P|\ge \deg(P)\ge 2\iota_X$ (\voir 4.1).
Il n'existe donc pas de collection primitive
d'ordre   $\le 2\iota_X-2$ dans $\fx$, de sorte que
$f_{j-1}=\binom{f_0}{j}$ pour tout $j \leq 2\iota_X-2$.  Puisque $X$ n'est
pas isomorphe
\`a
$ \PP ^{n}$, on a $2\iota_X-2< [n/2]+1
 $ par la proposition 4.11, ce qui contredit les hypoth\`eses.
\finpreuve

\medskip

\hs Nous sommes en mesure de prouver le

\medskip

\noindent {\bf Th\'eor\`eme \num} {\em
Soit $X$ une vari\'et\'e de Fano torique de dimension $n$,
de pseudo-indice $\iota_X$ et de nombre de Picard $\rho_X$.
Si $3 \iota_X\ge n+3$, on a $\rho_X(\iota_X-1)\le n$ (et donc $\rho_X\le
3$); de plus, si
$\rho_X(\iota_X-1)=n$, on a
$X\simeq (\PP ^{\iota_X-1})^{\rho_X}$.
}

\medskip

{\em D\'emonstration.} Le th\'eor\` eme se montre par r\'ecurrence
sur  $\rho_X$.
Si $\rho_X=1$, on a $X\simeq\PP ^{n}$ et le r\'esultat est d\'emontr\'e.
Sinon, d'apr\`es le lemme 4.6, il existe dans $X$ une relation
primitive contractible qui s'\'ecrit
$$ x_1+\dots+x_h=a_1 y_1+\dots +a_k y_k $$
avec $k\le \iota_X-2$. Soit
$A=V(\langle y_1,\dots,y_k\rangle)$ le lieu exceptionnel de la
contraction associ\'ee. D'apr\`es le lemme 4.5, $A$ est une
vari\'et\'e de Fano,
$\rho_A=\rho_X$  et $\iota_A\ge \iota_X-k$.
De plus, la contraction
fait de $A$ une fibration  $A\rightarrow B$ en $\PP ^{h-1}$,
o\` u $B$ est une vari\'et\'e de Fano torique satisfaisant
$\rho_B=\rho_X-1$, $\dim (B)=n-k-h+1$ et
$\iota_B\ge \iota_A\ge \iota_X-k$.\\
\hs V\'erifions l'in\'egalit\'e
$\iota_B>\frac{1}{3}\dim (B)+1$. Notons que $h\ge k+\iota_X\ge 2k+2$ et
que
$\dim (B)=n-k-h+1\le n-3k-1<n-3k$, d'o\`u
$$\iota_B\ge \iota_X-k\ge \frac{1}{3}n+1-k=\frac{1}{3}
(n-3k)+1>\frac{1}{3}\dim (B)+1.$$

\hs L'hypoth\`ese de r\'ecurrence donne donc $\rho_B(\iota_B-1)\le \dim
(B)$ et
$\rho_B=\rho_X-1\le 2$ puisque
$\iota_B>\frac{1}{3}\dim (B)+1$. Finalement:
$$(\rho_X-1)(\iota_X-k-1)\le \rho_B(\iota_B-1)\le\dim (B)=n-k-h+1$$
d'o\`u
$$\rho_X(\iota_X-1)-n\le \rho_X k+\iota_X-h-2k\le k+\iota_X-h\le 0.$$
\hs Supposons maintenant $\rho_X(\iota_X-1)=n$ et notons qu'il suffit de
montrer
 $A=X$, c'est-\`a-dire  $k=0$. En effet,
$X$ est alors une fibration en espaces projectifs sur une
vari\'et\'e de Fano torique
$B$ satisfaisant
$\rho_B \le 2$ et le corollaire 4.3 s'applique.
Supposons donc $k>0$. Les in\'egalit\'es pr\'ec\'edentes sont des \'egalit\'es,
d'o\`u
\[ \iota_B=\iota_A=\iota_X-k,\quad\rho_B(\iota_B-1)=\dim
(B),\quad\rho_X=\rho_A=3,\quad
\rho_B=2\quad\text{et}\quad h=k+\iota_X.\]
La vari\'et\'e $B$ est donc isomorphe \`a $\PP ^{\iota_X-k-1}\times\PP
^{\iota_X-k-1}$ (corollaire 4.3) et, puisque $A\rightarrow B$ est
une fibration en $\PP ^{\iota_X+k+1}$ et
$\iota_A=\iota_X-k$, la vari\'et\'e $A$  est isomorphe \`a
$\PP ^{\iota_X+k-1}\times\PP ^{\iota_X-k-1}\times\PP ^{\iota_X-k-1}$.
L'\'eventail $\Sigma_A$ est donc d\'etermin\'e par les relations
primitives
$$\overline{x}_1+\dots+\overline{x}_{\iota_X+k}=0,\
\overline{v}_1+\dots+\overline{v}_{\iota_X-k}=0\ \text{et}\
\overline{w}_1+\dots+\overline{w}_{\iota_X-k}=0,$$
o\`u l'on a not\'e $\overline{u}$
le g\'en\'erateur de $\Sigma_A$ induit par l'\'el\'ement
$u\in G(\fx)$ tel que $\langle u, y_1,\dots,y_k\rangle\in\fx$.
On a n\'ecessairement
$G(\fx)=\{x_1,\dots,x_{\iota_X+k},y_1,\dots,y_k,v_1,\dots,v_{\iota_X-k},w_1,
\dots,w_{\iota_X-k}\}$ puisque $\rho_X=3$.
Puisque toute classe extr\'emale de $X$ se restreint \`a une classe extr\'emale
dans $A$, il doit y avoir dans $X$ trois
relations primitives extr\'emales, de degr\'e au moins $\iota_{X}$, dont
les restrictions \`a $A$ sont les relations primitives ci-dessus.
La seule possibilit\'e est que les relations
$$ v_1+\dots+v_{\iota_X-k}+y_1+\dots+y_k=0 \mbox{ et }
w_1+\dots+w_{\iota_X-k}+y_1+\dots+y_k=0 $$
soient extr\'emales dans $X$, ce qui est absurde
si $k>0$ car ces relations ne sont pas disjointes
(\cite{contr}, Corollary 3.2).
\finpreuve

\subsection{Les petites dimensions.}
Dans ce paragraphe, on \'etudie les cas $\iota_X=[n/2]$ et $\iota_X=[n/2]-1$
lorsqu'ils
ne sont pas couverts par le th\'eor\`eme 4.7, c'est-\`a-dire
respectivement
pour
$n\le 7$ et $n\le 13$. Nous  en d\'eduirons que l'in\'egalit\'e
($\ast$) est
toujours verifi\'ee si $n\le 7$. Nous supposerons dans la suite $n\ge 4$.

\medskip

\noindent {\bf Proposition \num} {\em
Soit  $X$ une vari\'et\'e de Fano torique de
dimension $n\le13$ et
de pseudo-indice $\iota_X=[n/2]$ ou $\iota_X=[n/2]-1$.
Si $X$ n'est pas une fibration en $\PP ^{\iota_X-1}$,
l'entier
$\rho_X$ satisfait les in\'egalit\'es suivantes:

\smallskip

\renewcommand{\arraystretch}{1.9}
\scriptsize
\[
\begin{array}{|c|c|c|c|c|c|c|c|c|c|c|}
\cline{2-11}
\multicolumn{1}{c|}{}& n=4 & n=5 & n=6
& n=7 & n=8 & n=9 & n=10 & n=11 &n=12&n=13\\
\hline
\iota_X=[\frac{n}{2}] & \rho_X\le 2  & \rho_X\le 2
& \rho_X\le 2 & \rho_X\le 2
   \\
\hline
\iota_X=[\frac{n}{2}]-1 &
\multicolumn{2}{|c|}{}
 &\rho_X\le 4 & \rho_X\le 4 & \rho_X\le 3 &
\rho_X\le 3 & \rho_X\le 2  & \rho_X\le 3
& \rho_X\le 2 & \rho_X\le 2 \\
\cline{1-1}\cline{4-11}
\end{array}\]

\smallskip
}

\medskip

\hs V\'erifions que cette proposition entra\^ine le

\medskip

\noindent {\bf Corollaire \num} {\em
Soit  $X$ une vari\'et\'e de Fano torique de dimension
$n\le13$ et de pseudo-indice  $\iota_X=[n/2]$ ou $\iota_X=[n/2]-1$.
Alors $\rho_X(\iota_X-1)\le\dim (X)$ et
on a \'egalit\'e si et seulement si $X\simeq (\PP ^{\iota_X-1})^{\rho_X}$.
}

\medskip

{\em D\'emonstration.} L'in\'egalit\'e
$\rho_X(\iota_X-1)\le\dim (X)$, dans les cas non couverts par le
th\'eor\`eme 4.7, est \'equivalente aux bornes du tableau
ci-dessous:

\smallskip

{\renewcommand{\arraystretch}{1.9}
\scriptsize
\[
\begin{array}{|c|c|c|c|c|c|c|c|c|c|c|}
\cline{2-11}
\multicolumn{1}{c|}{}& n=4 & n=5 & n=6
& n=7 & n=8 & n=9 & n=10 & n=11 &n=12&n=13\\
\hline
\iota_X=[\frac{n}{2}] & \rho_X\le 4  &
\rho_X\le 5 & \rho_X\le 3 & \rho_X\le 3
   \\
\hline
\iota_X=[\frac{n}{2}]-1 &
\multicolumn{2}{|c|}{}
 &\rho_X\le 6 & \rho_X\le 7 & \rho_X\le 4 &
\rho_X\le 4 & \rho_X\le 3  & \rho_X\le 3
& \rho_X\le 3 & \rho_X\le 3\\
\cline{1-1}\cline{4-11}
\end{array}\]
}

\smallskip

\medskip

\hs Si $X$ n'est pas une fibration en $\PP ^{\iota_X-1}$, la
proposition 4.9 implique
$\rho_X(\iota_X-1)\le\dim (X)$ et  il n'y a jamais \'egalit\'e.\\
\hs Si $X$ est une fibration en $\PP ^{\iota_X-1}$ sur $Z$,
on a $\dim (Z)=n+1-\iota_X\le (n-1)/2$ et
$\iota_Z\ge \iota_X\ge \dim (Z)-1$. Si $\dim (Z) \ge 5$, la vari\'et\'e
$Z$ est un espace projectif (\voir 4.12) et donc $\rho_X=2$: le
r\'esultat est une cons\'equence du corollaire 4.3. Si $\dim
(Z)\le 4$, on a $\iota_Z\ge 3$ et $\rho_Z\le 2$, de sorte que le
corollaire 4.3 s'applique encore \`a $X$.
\finpreuve

\medskip

\hs Soit $P_X$ le polytope associ\'e \`a
$X$, c'est-\`a-dire le polytope simplicial convexe engendr\'e
par les \'el\'ements de $G(\fx)$. Les faces de $P_X$
correspondent aux c\^ones de $\fx$.
Soit $f_j$ le nombre de faces de dimension $j$ de $P_X$. Rappelons
l'\'egalit\'e
$f_0=|G(\fx)|=\rho_X+n$. L'\'eventail $\fx$ n'a pas de
collections primitives d'ordre $j<\iota_X$, autrement dit,
\begin{equation}\label{fj}
 f_{j-1}=\binom{f_0}{j} \ \   \mbox{ pour tout
}j<\iota_X.
\end{equation}
 De plus, si $\fx$ poss\`ede une collection primitive
$P=\{x_1,\dots,x_{\iota_X}\}$  d'ordre $\iota_X$, la
relation associ\'ee est n\'ecessairement
$x_1+\dots+x_{\iota_X}=0$,  puisque son degr\'e est au moins $\iota_X$.
Elle est contractible et la contraction associ\'ee
est une fibration en $\PP ^{\iota_X-1}$. Finalement:
$$ f_{\iota_X-1} < \binom{f_0}{\iota_X} \  \text{ si et seulement si } X
\text{ est une fibration en } \PP ^{\iota_X-1}.$$
\hs Nous renvoyons \`a \cite{mcmullenshephard} pour
les propri\'et\'es fondamentales des polytopes simpliciaux.
Rappelons en particulier
qu'il existe des relations lin\'eaires entre les $f_j$,
appel\'ees relations de Dehn-Sommerville,
de sorte que les nombres $f_0,\dots,f_{[n/2]-1}$ d\'eterminent tous
les $f_j$. Rappelons aussi le r\'esultat suivant (\cite{mcmullenshephard},
Chap. 2, Proposition 24).

\medskip

\noindent {\bf Proposition \num}
{\em Un polytope de dimension $n$ est un simplexe si et seulement si
$\displaystyle{f_{j-1}=\binom{f_0}{j}}
$ pour tout $ j\le [n/2]+1$.
}

\medskip

\marque Remarquons  en particulier  que si $\iota_X>[n/2]+1$, la
relation (\ref{fj}) et la proposition entra\^inent que
$X$ est isomorphe \`a $\PP ^{n}$.
Cela red\'emontre dans le cas torique
le th\'eor\`eme de Wi{\'s}niewski cit\'e dans l'introduction.

\medskip

\hs Les deux lemmes suivants donnent
des relations suppl\'ementaires entre les nombres $f_j$.

\medskip

\noindent {\bf Lemme \num} {\em
Soit $X$ une vari\'et\'e de Fano torique.
Si $\iota_X>1$, toutes les collections primitives de
$\fx$ d'ordre $\iota_X+1$ sont
deux \`a deux disjointes.
}

\medskip

{\em D\'emonstration.}
Soit $P=\{x_1,\dots,x_{\iota_X+1}\}$ une collection primitive d'ordre
$\iota_X+1$. Elle est n\'e\-ces\-sai\-re\-ment contractible puisque
$\deg( P)\le \iota_X+1 < 2\iota_X$ (\voir 4.1).
Ou bien $\deg (P)=\iota_X+1$, la relation primitive $r(P)$ est
$x_1+\dots+x_{\iota_X+1}=0$
et $P$ est  disjointe de
toutes les autres  collections primitives de
$\fx$ (\cite{contr}, Corollary 3.2), ou bien $\deg (P)=\iota_X$ et
la relation primitive
$r(P)$ est de la forme $x_1+\dots+x_{\iota_X+1}=y$.
Supposons par l'absurde qu'il existe deux
relations primitives
$$ x_1+\dots+x_r+y_1+\dots+y_s=z\ \ \mbox{ et }\ \
x_1+\dots+x_r+u_1+\dots+u_s=v $$
avec $r+s=\iota_X+1$ et $r>0$.
Comme $P$ est contractible, $\{y_1,\dots,y_s,v\}$ doit
contenir une collection primitive (\cite{contr}, Lemma 3.1), de sorte que
 $s+1\ge \iota_X $ et $r\le 2$. Si $r=2$ et $s=\iota_X-1$, la collection
$\{y_1,\dots,y_s,v\}$
est primitive, de relation  associ\'ee
$y_1+\dots+y_s+v=0$, ce qui est absurde. Si $r=1$ et
$s=\iota_X$, la collection $\{y_1,\dots,y_s,v\}$ est
primitive, de relation associ\'ee $y_1+\dots+y_s+v=w$. Par ce
qui pr\'ec\`ede, on a
$s=1$ et $\iota_X=1$, ce qui est exclu par l'hypoth\`ese.\finpreuve

\medskip

\noindent {\bf Lemme \num}  {\em
Si $X$ est une vari\'et\'e de Fano torique de dimension $n$, on a
$$ 12f_{n-3}\ge(3n+\iota_X-5)f_{n-2}.$$
}

{\em D\'emonstration.}
Posons $d_X=\sum\deg (N_{C/X})$, o\`u la somme porte sur toutes les
courbes invariantes de $X$. Si $C$ est une courbe invariante, on a $\deg
(N_{C/X})=-K_X\cdot C-2\ge \iota_X-2$, de sorte que $d_X\ge
(\iota_X-2)f_{n-2}$. D'autre part, on a $d_X=12f_{n-3}
-3(n-1)f_{n-2}$ (\cite{bat2}, Theorem 2.3.7), d'o\`u l'in\'egalit\'e
cherch\'ee.
\finpreuve

\medskip

\noindent {\em D\'emonstration de la proposition 4.9.}
Posons $k=[n/2]$ et supposons $\iota_X=k$.
Comme $X$ n'est pas une fibration en $\PP ^{\iota_X-1}$, on a
$f_{j-1}=\binom{f_0}{j}$ pout tout $j\le k$. Les relations de Dehn-Sommerville
permettent d'exprimer  $f_k$ \`a l'aide
de $f_0$
(\cite{mcmullenshephard}, \S\ 2.4):
\begin{align*}
f_k&=\sum_{j=0}^{k-1}(-1)^{k-j-1}
\frac{j+1}{k+1}\binom{2k-j}{k}\binom{f_0}{j+1}\quad\text{ si }n
\text{ est pair,} \\
 f_k&=\sum_{j=-1}^{k-1}(-1)^{k-j-1}\binom{2k-j+1}{k+1}\binom{f_0}{j+1}
\quad\text{ si }n
\text{ est impair.}
\end{align*}
Le lemme 4.13 donne l'in\'egalit\'e
\[ \binom{f_0}{k+1}-f_{k}\le \frac{f_0}{k+1},  \]
d'o\`u les bornes de la premi\`ere ligne du tableau.\\
\hs Supposons $\iota_X=k-1$ et que $X$ n'est pas
une fibration en $\PP^{\iota_X-1}$. On a alors
$$ f_{j-1}=\binom{f_0}{j} \text{ pour tout }j\le k-1\quad\text{et}\quad
\binom{f_0}{k}-f_{k-1}\le \frac{f_0}{k}.$$
Les relations de Dehn-Sommerville permettent d'exprimer
tous les $f_j$ avec $j\ge k$ \`a l'aide
de $f_0$ et $f_{k-1}$. Si $n$ est pair,
on a (\cite{mcmullenshephard}, \S\ 2.4)
\begin{align*}
f_{n-2}&=kf_{k-1}+
\sum_{j=1}^{k-1}(-1)^j\frac{k-j}{k+j-1}\left((k-1)\binom{k+j}{k}+
\binom{k+j-1}{k}\right)
\binom{f_0}{k-j}\quad\text{et}\\
f_{n-3}&=\binom{k}{2}f_{k-1}+\sum_{j=1}^{k-1}(-1)^j\frac{k-j}{k+j-2}
\left(\binom{k}{2}\binom{k+j}{k}\right.\\
& \hspace{6.5cm}\left.
+(k-2)\binom{k+j-1}{k}+\binom{k+j-2}{k}\right)\binom{f_0}{k-j}.
\end{align*}
Et si $n$ est impair, on a
\begin{align*}
f_{n-2}&=(2k+1)f_{k-1}+\sum_{j=1}^k(-1)^j\frac{2k+1}{k+j}
\left(k\binom{k+j+1}{k+1}+\binom{k+j}{k+1}\right)
\binom{f_0}{k-j}\quad\text{et} \\
f_{n-3}&=k^2f_{k-1}
+\sum_{j=1}^k(-1)^j\frac{2k}{k+j-1}\left(\binom{k}{2}
\binom{k+j+1}{k+1}\right.\\
& \hspace{7cm}\left.
+(k-1)\binom{k+j}{k+1}+\binom{k+j-1}{k+1}\right)\binom{f_0}{k-j}.
\end{align*}
Finalement, en appliquant le lemme 4.14, on obtient
$f_{k-1}\le\Psi_{k}(f_0)$, o\`u, si $n$ est pair,
\begin{multline*}
\Psi_{k}(f_0)=\frac{k^3-k^2+12}{k^2}\binom{f_0}{k-1}\\
+\sum_{j=2}^{k-1}\frac{(-1)^{j-1}(k-j)(k+j-3)!}{k(k!)(j!)}
\Bigl((k+j-1)(k+j-2)k+12j
\Bigr) \binom{f_0}{k-j}
\end{multline*}
et, si $n$ est impair,
\begin{multline*}
\Psi_{k}(f_0)=\frac{2k^3+3k^2-2k+21}{(k-1)(2k+3)}\binom{f_0}{k-1}
+\sum_{j=2}^{k}\frac{(-1)^{j-1}(k+j-2)!}{(k-1)(2k+3)(k!)(j!)}\Bigl(2k^4+(4j-1)k^
3\\
+2(j^2-2)k^2+(j^2+17j+3)k+3j(1-j)\Bigr)\binom{f_0}{k-j}.
\end{multline*}

On d\'eduit alors du lemme 4.13 l'in\'egalit\'e
$$ \binom{f_0}{k}-\frac{f_0}{k}-\Psi_k(f_0)\le 0,$$
ce qui donne les bornes de la deuxi\`eme ligne du tableau par
une \'etude directe.
\finpreuve

\section{Cha\^ines de courbes rationnelles}

\hs Soit $X$ une vari\'et\'e projective, lisse et connexe et
soient $V^1,\dots,V^k$ des familles propres irr\'eductibles de courbes
rationnelles irr\'eductibles sur $X$ (\voir 3.14).
Fixons $x \in \lieu (V^1)$, posons $\lieu (V^1)_x =\lieu ( V^1_x)$ et
$Z_1 = {\mathcal U}^1_x$ (\voir 3.14) et, avec les notations
\begin{equation*}
\begin{CD}
\mathcal U^k @){\ev_k})) X \\
@V{\pi_k}VV \\
V^k
\end{CD}
\end{equation*}
d\'efinissons par r\'ecurrence pour
$k\ge2$:
\begin{eqnarray*}
Z_k &=&  \pi_k ^{-1}
(\pi_k (\ev_k ^{-1} ( \lieu  (V^1,\dots,V^{k-1})_x \cap \lieu (V^k) )))\\
\lieu  (V^1,\dots,V^k)_x &=&  \ev_k ( Z_k)\\
U^k_x&=&\lieu (V^1,\dots,V^k)_x \setminus
\lieu (V^1,\dots,V^{k-1})_x
\end{eqnarray*}
Le ferm\'e $\lieu  (V^1,\dots,V^k)_x$
est donc l'ensemble des points $y$ de $X$ tels qu'il
existe des courbes $C^1,\dots,C^k$ avec
\begin{itemize}
\item[$\bullet$]$[C^j] \in V^j$;
\item[$\bullet$] les intersections $C^1 \cap C^2,\dots, C^{k-1} \cap C^k $
ne sont pas vides;
\item[$\bullet$]
$x\in C^1$ et $y\in C^k$.
\end{itemize}

\medskip

\noindent {\bf Lemme \num} {\em Si les classes
des familles $V^1,\dots,V^k$ dans $\Nun_{1}(X)_{\QQ}$
sont lin\'eairement ind\'ependantes et que  $\lieu
(V^1,\dots,V^k)_x$ n'est pas vide,
\begin{enumerate}
\item[(a)] le morphisme d'\'evaluation $\ev_k : Z_k \to \lieu
(V^1,\dots,V^k)_x $ est fini au-dessus de l'ouvert   $U^k_x$, qui n'est
pas vide;
\item[(b)]  toute courbe trac\'ee
sur $\lieu (V^1,\dots,V^k)_x $ est alg\'ebriquement
\'equivalente dans  cette vari\'et\'e \`a une combinaison
lin\'eaire
\`a coefficients rationnels de
courbes dans $
V^1,\dots,V^k$.
\end{enumerate}
}

\medskip

{\em D\'emonstration.} Le r\'esultat se montre par r\'ecurrence sur l'entier
$k\ge 1$. Si $k=1$, le point (a) est une cons\'equence du
lemme de cassage de Mori
(\cite{cass}, Theorem 6) et le
point (b) est le lemme 3.15. Supposons $k \ge 2$.

\hs D\'emontrons (a).  Soit $z $ un point de $
\lieu (V^1,\dots,V^{k-1})_x\cap \lieu (V^k)$, qui n'est pas vide par
hypoth\`ese. Si $\lieu (V^k_z)\subset\lieu (V^1,\dots,V^{k-1})_x$, les
classes des familles $V^1,\dots,V^k$ sont li\-n\'eai\-rement
d\'e\-pen\-dantes dans
$\Nun_{1}(X)_{\QQ}$ par
hypoth\`ese de r\'ecurrence, ce qui est absurde.
Ainsi, $U^k_x$ n'est pas vide.
Si un point $y $ de $ U^k_x$ v\'erifie
$\dim (\ev_k^{-1}(y) \cap Z_k) \ge 1$,
\begin{itemize}
\item[$\bullet$]soit  $\dim (\lieu (V^k_y) \cap
\lieu (V^1,\dots,V^{k-1})_x ) \ge 1$ et, \`a nouveau par hypoth\`ese de
r\'ecurrence, les classes des familles $V^1,\dots,V^k$ sont
lin\'eairement d\'ependantes dans $\Nun_{1}(X)_{\QQ}$, ce qui contredit
l'hypoth\`ese;
\item[$\bullet$]soit il existe une famille de
dimension 1 de courbes rationnelles de $V^{k}$ passant toutes par deux
points distincts de $X$ fix\'es, ce qui, par le lemme de cassage
(\cite{cass}, Theorem 6), est absurde puisque la famille
$V^{k}$ est propre.
\end{itemize}

\hs D\'emontrons   (b). Soit $C$ une
courbe trac\'ee sur $\lieu (V^1,\dots,V^k)_x$ telle que
$[C] \notin V^k$. Si $C$ est dans
$\lieu (V^1,\dots,V^{k-1})_x$, l'hypoth\`ese de r\'ecurrence
permet de conclure. Supposons donc que $C$
n'est pas contenue dans $\lieu (V^1,\dots,V^{k-1})_x$.
 Soit ${\mathcal C }\subset \ev _k^{-1}(C) \cap Z_k$ une
courbe irr\'eductible
dominant $C$, soit ${\mathcal S}\subset {\mathcal U }^k$ la surface
irr\'eductible
$\pi_k ^{-1} (\pi_k ({\mathcal C}))$, soit $S$ la surface $\ev _k
({\mathcal S})
\subset X$ et soit ${\mathcal C}'\subset {\mathcal S}$ une courbe
dominant $ S \cap
\lieu (V^1,\dots,V^{k-1})_x$.
Toute courbe trac\'ee sur
${\mathcal S}$ est alg\'ebriquement
\'equivalente dans ${\mathcal S}$ \`a une   combinaison lin\'eaire
\`a coefficients rationnels  de la multisection ${\mathcal C}'$
et d'une fibre de ${\pi_k}_{\vert{\mathcal S}} : {\mathcal S} \to \pi_k
({\mathcal S})$  (\cite{Kol99}, II.4.19).
Toute courbe trac\'ee
sur $S$ est donc alg\'ebriquement
\'equivalente dans $S$, donc dans $\lieu (V^1,\dots,V^k)_x $, \`a une
combinaison lin\'eaire
\`a coefficients rationnels  de $\ev_k({\mathcal C}')$ et d'une
courbe de $V^k$
(\cite{Kol99}, II.4.4.2).
 Il
reste
\`a remarquer que
$\ev_k({\mathcal C}') \subset \lieu (V^1,\dots,V^{k-1})_x$:
l'hypoth\`ese de r\'ecurrence
permet de conclure.
\finpreuve

\medskip

\hs On en d\'eduit le r\'esultat suivant.

\medskip

\noindent {\bf Th\'eor\`eme \num} {\em
Si les classes des familles $V^1,\dots,V^k$ dans $\Nun_{1}(X)_{\QQ}$
sont lin\'eairement ind\'e\-pendantes,
 $\lieu (V^1,\dots,V^k)_x$ est vide ou de dimension $\ge
 - \sum_{j=1}^k K_X \cdot V^j
- k $.
}

\medskip

{\em D\'emonstration.} On proc\`ede par r\'ecurrence sur l'entier
$k\ge 1$. Si $k=1$, l'estimation sur la dimension a d\'eja \'et\'e
mentionn\'ee en 3.15. Supposons $k\ge 2$ et que $\lieu
(V^1,\dots,V^k)_x$  n'est pas vide.
 Le lemme
pr\'ec\'edent donne
$$ \dim (\lieu (V^1,\dots,V^k)_x)\ge \dim (Z_k).$$
Si $y$ est un point g\'en\'eral de $\lieu  (V^1,\dots,V^{k-1})_x \cap
\lieu (V^k) $, on a
$\dim (\ev_k^{-1}(y)) = \dim (V_y^k)$ et, si l'on note
$$W_k = \ev_k^{-1}(\lieu  (V^1,\dots,V^{k-1})_x \cap \lieu (V^k)),$$
on a, en utilisant l'hypoth\`ese de r\'ecurrence,
\begin{eqnarray*}
\dim (W_k)
&=& \dim (V_y^k) +\dim(\lieu  (V^1,\dots,V^{k-1})_x \cap \lieu
(V^k))\\
&\ge& \dim (V_y^k) - \sum_{j=1}^{k-1} K_X \cdot V^j -(k-1) + \dim
(\lieu (V^k) )-n\\
 &\ge& \dim (V_k) - \sum_{j=1}^{k-1} K_X \cdot V^j -k +2 -n,
\end{eqnarray*}
d'o\`u, par l'in\'egalit\'e (\ref{dim})  de 3.14,
$$ \dim (W_k) \ge - \sum_{j=1}^{k} K_X \cdot V^j -k -1.$$
Comme $\dim (Z_k)=\dim (W_k) + 1 $, le lemme pr\'ec\'edent
permet de conclure.
\finpreuve

\medskip

\hs On d\'eduit du th\'eor\`eme le r\'esultat suivant.

\medskip

\noindent {\bf Corollaire \num} {\em
Soit $X$ une vari\'et\'e de Fano de nombre de Picard $\rho_X$ et
de pseudo-indice $\iota_X$. S'il existe des familles propres
irr\'eductibles $V^1,\dots,V^{\rho_X}$ de cour\-bes rationnelles
irr\'e\-ductibles sur $X$ dont les classes dans
$\Nun_{1}(X)_{\QQ}$ sont lin\'eairement ind\'ependantes et que
$\lieu (V^1,\dots,V^{\rho_X})_x $ n'est pas vide, on a
$\rho_X (\iota_X -1)\le \dim(X)$.
}

\medskip

\hs Il n'est \'evidemment pas facile d'assurer l'existence de
familles propres de courbes rationnelles irr\'eductibles sur $X$
v\'erifiant les conditions du corollaire pr\'ec\'edent.
Si $R\subset\NE (X)$ est une ar\^ete,
les courbes rationnelles irr\'eductibles dont la classe
appartient \`a $R$ et de degr\'e anticanonique minimal forment
une famille propre. En consid\'erant ces familles de
courbes rationnelles, on montre le

\medskip

\noindent {\bf Corollaire \num} {\em
Soit $X$ une vari\'et\'e de Fano homog\`ene. On a $\rho_X (\iota_X -1) \le
\dim(X)$. }

-----------

\noindent L.B. {\em e-mail: bonavero@ujf-grenoble.fr
}

\noindent S.D. {\em e-mail: druel@ujf-grenoble.fr
}

\noindent {\em Institut Fourier, UFR de Math\'ematiques,
Universit\'e de Grenoble 1,
UMR 5582,
BP 74,
38402 Saint Martin d'H\`eres,
FRANCE.

}

\noindent C.C. {\em e-mail: ccasagra@mat.uniroma1.it
}

\noindent {\em Universit\`a di Roma ``La Sapienza'',
Dipartimento di Matematica, Piazzale Aldo Moro, 2,
00185 Rome, ITALIE.

}

\noindent O.D. {\em e-mail: debarre@math.u-strasbg.fr
}

\noindent {\em Math\'ematique - IRMA - UMR 7501,
Universit\'e Louis Pasteur,
7, rue Ren\'e Descartes,
67084 Strasbourg Cedex,
FRANCE.

}


\medskip


\begin{thebibliography}{FATA268}


\bibitem[An85]{And85}T.\ Ando. On extremal rays of the higher dimensional
varieties.
{\em Invent. Math. {\bf81} (1985), 347--357.}

\bibitem[ABW92]{ABW92}M.\ Andreatta, E.\ Ballico, J.A.\ Wi\'sniewski.
Vector bundles and adjunction.
{\em Internat. J. Math. {\bf3} (1992), 331--340.}

\bibitem[AW98]{AW98} M.\ Andreatta, J.A.\ Wi\'sniewski.
Contractions of smooth varieties. II. Computations and
applications.
{\em Boll. Unione Mat. Ital. Sez. B Artic. Ric. Mat. (8) 1 (1998),
 343--360.}

\bibitem[Ba91]{bat1} V.V.\ Batyrev.
On the classification of smooth projective toric varieties.
{\em Tohoku Mathematical Journal {\bf43} (1991), 569--585.}

\bibitem[Ba99]{bat2}
V.V.\ Batyrev.
On the classification of toric Fano $4$-folds.
{\em Journal of Mathematical Sciences (New York) {\bf94} (1999),
1021--1050.}

\bibitem[BCW01]{BCW01}L.\ Bonavero, F.\ Campana, J.A.\ Wi\'sniewski.
Vari\'et\'es projectives complexes dont l'\'eclat\'ee en
un point est de Fano.
{\em C.R. Acad. Sci. Paris, Ser. I {\bf334} (2002), 463--468.}

\bibitem[Ca01]{contr}
C.\ Casagrande.
Contractible classes in toric varieties.
{\em Pr\'e\-pu\-blica\-tion math.AG/0111332 (2001).}

\bibitem[CMS00]{CMS00}K.\ Cho,
Y.\ Miyaoka, N.\ Shepherd-Barron.
Characterizations of projective spaces and applications.
{\em Preprint (2000).}

\bibitem[Fu93]{Ful93}W.\ Fulton.
Introduction to toric varieties.
{\em Annals of mathematics studies {\bf131},
Princeton University Press, 1993.}



\bibitem[Gr68]{gro}A. Grothendieck.
Le groupe de Brauer. I. Alg\`ebres d'Azumaya et interpr\'etations
diverses. {\em Dix expos\'es sur la cohomologie des sch\'emas,
46--66, North-Holland,  Amsterdam; Masson, Paris, 1968. }



\bibitem[IP99]{IsP99}
V.\ A.\ Iskovskikh, Yu.\ G.\ Prokhorov.
Fano varieties.
Algebraic geometry, V.
{\em  Encyclopaedia Math. Sci. {\bf47}, Springer-Verlag, Berlin,
1999.}

\bibitem[Ka89]{Kaw89}Y.\ Kawamata. Small contractions of four
dimensional algebraic manifolds.
{\em Math. Ann. {\bf284} (1989), 595--600}.

\bibitem[Ke01]{Keb01}S.\ Kebekus.
Characterizing the projective space after Cho, Miyaoka and
Shepherd-Barron.
{\em Pr\'e\-pu\-blication math.AG/0107069. \`A para\^itre dans \og
Festschrift in honor of Hans Grauert\fg.}


\bibitem[Kl88]{Kle88}P.\ Kleinschmidt.
A classification of toric varieties with few generators.
{\em Aequationes Math. {\bf35} (1988),  254--266.}

\bibitem[Ko96]{Kol99}J.\ Koll\'ar.
Rational curves on algebraic varieties.
{\em Ergebnisse der Mathematik und ihre Grenzgebiete. {\bf3} Folge
{\bf032},
 Springer-Verlag, 1996}.


\bibitem[KMM92]{kmm}   J.  Koll\'ar, Y.  Miyaoka, S.
Mori.  Rational Connectedness and Boundedness of  Fano
Manifolds. {\em J. Diff. Geom. {\bf 36} (1992), 765--769.}

\bibitem[La83]{Laz83}
R.\ Lazarsfeld. Some applications of the
theory of positive vector bundles.
{\em Complete intersections (Acireale,
1983), 29-61, Lecture Notes in Math. {\bf1092}, Springer-Verlag, Berlin,
1984.}


\bibitem[MMS71]{mcmullenshephard}
P.\ McMullen, G.C.\ Shephard.
Convex Polytopes and the Upper Bound Conjecture.
{\em  London Mathematical Society Lecture Note Series {\bf3},
Cambridge University Press, 1971.}

\bibitem[Mo79]{cass}  Mori, S., Projective manifolds with ample
tangent bundles, {\it Ann. of Math.} {\bf 110} (1979), 593--606.

\bibitem[Mo82]{Mor82}S.\ Mori. Threefolds whose canonical bundles are
not numerically effective.
{\em Ann. of Math. {\bf116} (1982), 133--176.}

\bibitem[MM81]{MM81}
S.\ Mori, S.\ Mukai.
Classification of Fano $3$-folds with $B\sb{2}\ge 2$.
{\em Manuscripta Math. {\bf36}
(1981), 147-162.}


\bibitem[Mu88]{Muk88} S.\ Mukai.
Birational geometry of algebraic varieties. Open problems.
{\em The 23rd Intern. Sympos. Division of Math. The Taniguchi
Foundation. Kataka (1988).}

\bibitem[Od88]{Oda88}T.\ Oda. Convex bodies and
algebraic geometry: an introduction to the theory
of toric varieties.
{\em
Ergebnisse der Mathematik und ihrer Grenzgebiete. 3. Folge. {\bf015},
Springer-Verlag, 1988.}

\bibitem[Wi90]{wisn2}
J.A.\ Wi{\'s}niewski.
On a conjecture of Mukai.
{\em Manuscripta Mathematica {\bf68} (1990), 135--141.}


\bibitem[Wi91]{Wis91}
J.A.\ Wi\'sniewski. On contractions of extremal rays
of Fano manifolds.
{\em J. Reine Angew. Math. {\bf417}
(1991), 141--157.}


\end{thebibliography}
\end{document}